\makeatletter
\IfFileExists{./numapde-latex/numapde-packages.sty}{\providecommand*{\input@path}{}\edef\input@path{{./numapde-latex/}\input@path}}{}
\makeatother

\RequirePackage{algpseudocode}
\documentclass{numapde-preprint}

\usepackage{numapde-semantic}
\usepackage{numapde-style}
\usepackage{numapde-local}

\IfFileExists{./numapde-bibliography/numapde.bib}{\addbibresource{./numapde-bibliography/numapde.bib}}{\addbibresource{numapde.bib}}
\hypersetup{
	pdftitle={Optimal Control of the Kirchhoff Equation},
	pdfauthor={Masoumeh Hashemi, Roland Herzog, Thomas M. Surowiec},
	pdfkeywords={PDE-constrained optimization, optimal control, nonlocal equation, Kirchhoff equation, quasilinear equation, semismooth Newton method}
}

\title{Optimal Control of the Kirchhoff Equation\thanks{This work was supported by a DFG grant HE 6077/8-1 within the Priority Program SPP 1962 (Non-smooth and Complementarity-based Distributed Parameter Systems: Simulation and Hierarchical Optimization), which is gratefully acknowledged.}}
\subtitle{}
\shorttitle{Optimal Control of the Kirchhoff Equation}

\author{Masoumeh Hashemi\thanks{Interdisciplinary Center for Scientific Computing, Heidelberg University, 69120 Heidelberg, Germany (\email{masoumeh.hashemi@iwr.uni-heidelberg.de}, \url{https://scoop.iwr.uni-heidelberg.de/team/masoumeh-hashemi}, \orcid{0000-0002-2835-2249}, \email{roland.herzog@iwr.uni-heidelberg.de}, \url{https://scoop.iwr.uni-heidelberg.de/team/roland-herzog}, \orcid{0000-0003-2164-6575}).}
\and
Roland Herzog\footnotemark[2]
\and
Thomas M. Surowiec\thanks{Philipps University Marburg, Fachbereich Mathematik und Informatik, Hans-Meerwein-Straße 6, 35032 Marburg, Germany (\email{surowiec@mathematik.uni-marburg.de}, \url{https://www.mathematik.uni-marburg.de/\textasciitilde surowiec}, \orcid{0000-0003-2473-4984}).}}
\shortauthor{M. Hashemi, R. Herzog and T. M. Surowiec}

\dedication{}

\IfFileExists{numapde-uncolorizeHyperrefIfChanges.sty}{\RequirePackage{numapde-uncolorizeHyperrefIfChanges}}{}

\begin{document}
\maketitle

\begin{abstract}
We consider an optimal control problem for the steady-state Kirchhoff equation, a prototype for nonlocal partial differential equations, different from fractional powers of closed operators.
Existence and uniqueness of solutions of the state equation, existence of global optimal solutions, differentiability of the control-to-state map and first-order necessary optimality conditions are established.
The aforementioned results require the controls to be functions in $H^1$ and subject to pointwise upper and lower bounds.
In order to obtain the Newton differentiability of the optimality conditions, we employ a Moreau-Yosida-type penalty approach to treat the control constraints and study its convergence.
The first-order optimality conditions of the regularized problems are shown to be Newton diffentiable, and a generalized Newton method is detailed.
A discretization of the optimal control problem by piecewise linear finite elements is proposed and numerical results are presented.\end{abstract}

\begin{keywords}
PDE-constrained optimization, optimal control, nonlocal equation, Kirchhoff equation, quasilinear equation, semismooth Newton method\end{keywords}

\begin{AMS}
\href{https://mathscinet.ams.org/msc/msc2010.html?t=49J20}{49J20}, \href{https://mathscinet.ams.org/msc/msc2010.html?t=49K20}{49K20}, \href{https://mathscinet.ams.org/msc/msc2010.html?t=35J62}{35J62}, \href{https://mathscinet.ams.org/msc/msc2010.html?t=47J05}{47J05}
\end{AMS}

\section{Introduction}
\label{section:introduction}
In this paper we study an optimal control problem governed by a nonlinear, nonlocal partial differential equation (PDE) of Kirchhoff-type
\begin{equation}\label{eq:Kirchhoff_equation_general}
	\left\{
		\begin{aligned}
			-M\paren[auto](){x, \norm{\nabla y}_{L^2(\Omega)}^2; u} \laplace y & = f & & \text{in } \Omega,
			\\
			y & = 0 & & \text{on } \partial \Omega.
		\end{aligned}
	\right.
\end{equation}
Here, $\Omega \subset \R^N$ is an open and bounded set and the right-hand side $f$ belongs to $L^2(\Omega)$.
We focus on the particular case $M\paren[auto](){x,s;u} = u(x) + b(x) \, s$, which has been considered previously, \eg, in \cite{FigueiredoMoralesRodrigoSantosSuarez:2014:1,DelgadoFigueiredoGayteMoralesRodrigo:2017:1}.
Here $u$ and $b$ are strictly positive functions and $u$ serves as the control.
The full set of assumptions is given in \cref{section:existence_of_an_optimal_solution}.
We mention that in case $u$ and $b$ are positive \emph{constants}, \eqref{eq:Kirchhoff_equation_general} has a variational structure; see \cite{FigueiredoMoralesRodrigoSantosSuarez:2014:1}.

Equation \eqref{eq:Kirchhoff_equation_general} is the steady-state problem associated with its time-dependent variant
\begin{equation}\label{eq:Kirchhoff_equation_time-dependent}
	\paren[auto]\{.{%
		\begin{aligned}
			y_{tt} - M\paren[auto](){x, \norm{\nabla y}_{L^2(\Omega)}^2; u} \laplace y & = f & & \text{in } \Omega \times (0,T),
			\\
			y & = 0 & & \text{on } \partial \Omega \times (0,T),
			\\
			y(x,0) = y_0(x), \quad y_t(x,0) & = y_1(x) & & \text{in } \Omega.
		\end{aligned}
	}
\end{equation}
Problem \eqref{eq:Kirchhoff_equation_time-dependent} models small vertical vibrations of an elastic string with fixed ends, when the density of the material is not constant. 
Specifically, the control~$u$ is proportional to the inverse of the string's cross section; see \cite{Ma:2005:1,FigueiredoMoralesRodrigoSantosSuarez:2014:1}.

PDEs with nonlocal terms play an important role in physics and technology and they can be mathematically challenging. 
Although in some cases variational reformulations are available, the models \eqref{eq:Kirchhoff_equation_general}, \eqref{eq:Kirchhoff_equation_time-dependent} do not allow this in general. 
Thus, despite the deceptively simple structure, \eqref{eq:Kirchhoff_equation_general} requires a set of analytical tools not often employed in PDE-constrained optimization.
Existence and uniqueness of solutions for \eqref{eq:Kirchhoff_equation_general} have been investigated in \cite{FigueiredoMoralesRodrigoSantosSuarez:2014:1} and \cite{DelgadoFigueiredoGayteMoralesRodrigo:2017:1}; see also the references therein.
For further applications of nonlocal PDEs, we refer the reader to \cite{Eringen:1983:1,AhmedElgazzar:2007:1,KavallarisSuzuki:2018:1}.

\cite{DelgadoFigueiredoGayteMoralesRodrigo:2017:1} studied an optimal control problem for \eqref{eq:Kirchhoff_equation_general} with the following cost functional 
\begin{equation}
	\label{eq:cost_functional_Delgado_et_al}
	J(y,u) 
	= 
	\frac{1}{2} \norm{y-y_d}_{L^2(\Omega)}^2 + \frac{\lambda}{2} \norm{u}_{L^2(\Omega)}^2 
\end{equation}
with an admissible set $\Uad = \setDef{u \in L^2(\Omega)}{u \ge u_a > 0 \text{ \ale in } \Omega}$. 
However we believe that the proof of existence of an optimal solution in this work has a flaw.
We give further details in the appendix.
Moreover, the proof in \cite{DelgadoFigueiredoGayteMoralesRodrigo:2017:1} is explicitly tailored to such tracking type functionals.
In the present work we see it necessary to modify the control cost term to contain the stronger $H^1$-norm.
We also allow for a more general state dependent term, which leads to the objective
\begin{equation}
	\label{eq:cost_functional}
	J(y,u) 
	= 
	\int_\Omega \varphi(x,y(x)) \dx + \frac{\lambda}{2}\norm{u}_{H^1(\Omega)}^2
\end{equation}
and a set of admissible controls in $H^1(\Omega)$.
In this setting, we prove the weak-strong continuity of the control-to-state operator into $H^1_0(\Omega) \cap W^{2,p}(\Omega)$ for any $p \in [1,\infty)$. 
Moreover, we work with a pointwise lower bound on admissible controls.
This bound has an immediate technological interpretation, representing an upper bound on the string's cross section.
On the other hand, in order to prove the existence of a globally optimal solution, we impose an additional upper bound on the admissible controls.
We are able to prove the Fréchet differentiability of the control-to-state map so that we can derive optimality conditions in a more straightforward way than by the Dubovitskii-Milyoutin formalism utilized in \cite{DelgadoFigueiredoGayteMoralesRodrigo:2017:1}.

The first-order optimality conditions obtained when minimizing \eqref{eq:cost_functional} subject to \eqref{eq:Kirchhoff_equation_general} results in a variational inequality of nonlinear obstacle type in $H^1$.
Unfortunately, it is not known whether a formulation of this condition exists which is differentiable in a generalized sense so that a generalized Newton method can be applied.
Therefore, we choose to relax and penalize the bound constraints via a Moreau-Yosida regularization, which amounts to a quadratic penalty of the bound constraints for the control.
In this setting, we can prove the generalized (Newton) differentiability of the optimality system.
A similar philosophy, albeit for a different problem, has been pursued by \cite{AdamHintermuellerSurowiec:2018:1}.
We also mention \cite[Chapter~9.2]{Ulbrich:2011:1} for an approach via a regularized dual obstacle problem.
Relaxing the lower and upper bounds, however, adds new difficulties, since the existence of a solution of the Kirchhoff equation \eqref{eq:Kirchhoff_equation_general} can only be guaranteed for positive controls.
Therefore, we compose the control-to-state map with a smooth cut-off function.
We then study the convergence of global minimizers as the penalty parameter goes to zero, see \cref{theorem:optimality_system_penalized_problem} for details.
We can expect a corresponding result to hold also for locally optimal solutions under an assumption of second-order sufficient optimality conditions, but this is not investigated here.

To summarize our contributions in comparison to \cite{DelgadoFigueiredoGayteMoralesRodrigo:2017:1}, we consider a more general objective, present a simpler proof for the existence of a globally optimal control, prove the differentiability of the control-to-state map and generalized differentiability of the optimality system for a regularized version of the problem as well as the applicability of a generalized Newton scheme.
We also describe a structure preserving finite element discretization of the problem and the discrete counterpart of the generalized Newton method.

The paper is organized as follows. 
In \cref{section:existence_of_an_optimal_solution}, we review existence and uniqueness results for solutions of the Kirchhoff equation \eqref{eq:Kirchhoff_equation_general} and prove the existence of a globally optimal control.
Subsequently, we prove the Fréchet differentiability of control-to-state operator and derive a system of necessary optimality conditions for a regularized problem in \cref{section:optimality_system}. 
In \cref{section:generalized_Newton_method}, we prove the Newton differentiability of the optimality system and devise a locally superlinealy convergent scheme in appropriate function spaces.
\Cref{section:discretization} addresses the discretization of the optimal control problem, its optimality system and the generalized Newton method by a finite element scheme.
The paper concludes with numerical results in \cref{section:numerical_experiments}.

\section{Optimal Control Problem: Existence of a Solution}
\label{section:existence_of_an_optimal_solution}

In this work we are interested in the study of the following optimal control problem for a stationary nonlinear, nonlocal Kirchhoff equation:
\begin{subequations}
	\label{eq:optimal_control_problem}
	\begin{align}
		\label{eq:optimal_control_problem_objective}
		\text{Minimize} \quad 
		&
		J(y,u)
		\coloneqq
		\int_\Omega \varphi(x,y(x)) \dx 
		+
		\frac{\lambda}{2} \norm{u}_{H^1(\Omega)}^2
		\\
		\label{eq:optimal_control_problem_Kirchhoff_equation}
		\text{subject to} \quad
		&
		\paren[auto]\{.{%
			\begin{aligned}
				-\paren[auto](){u + b \, \norm{\nabla y}_{L^2(\Omega)}^2}\laplace y & = f & & \text{in } \Omega, 
				\\
				y & = 0 & & \text{on } \partial\Omega
			\end{aligned}
		}
		\\
		\label{eq:optimal_control_problem_admissible_set}
		\text{and} \quad
		&
		u \in \Uad \cap H^1(\Omega)
		.
	\end{align}
\end{subequations}
The set of admissible controls is given by
\begin{equation}
	\label{eq:admissible_set}
	\Uad 
	= 
	\setDef{u \in L^\infty(\Omega)}{u_a(x) \le u(x) \le u_b(x) \text{ \ale\ in } \Omega}
	.
\end{equation}

The following are our standing assumptions.
\begin{assumption}\label{assumption:general}
	We assume that $\Omega \subset \R^N$ is a bounded domain of class~$C^{1,1}$ with $1 \le N \le 3$; see for instance \cite[Chapter~2.2.2]{Troeltzsch:2010:1}
	The right-hand side~$f$ is a given function in $L^\infty(\Omega)$ satisfying $f \ge f_0$ \ale, where $f_0$ is a positive real number.
	The bounds $u_a$ and $u_b$ are functions in $C(\overline \Omega)$ such that $u_b \ge u_a \ge u_0$ holds for some positive real number $u_0$.
	Finally, we assume $b \in W^{1,\infty}(\Omega)$ with $b \ge b_0$ \ale for some positive real number $b_0$.
\end{assumption}

The integrand~$\varphi$ in the objective is assumed to satisfy the following standard assumptions; see for instance \cite[Chapter~4.3]{Troeltzsch:2010:1}:
\begin{assumption}\label{assumption:Nemytskii}
	\begin{enumeratearabic}
		\item 
			$\varphi \colon \Omega\times \R\rightarrow\R$ is Carathéodory and of class~$C^2$, \ie,
			\begin{enumerate}
				\item 
					$\varphi(\cdot, y) \colon \Omega\rightarrow \varphi(x,y)$ is measurable for all $y\in \R$,

				\item 
					$\varphi(x,\cdot) \colon \R \rightarrow \varphi(x,y)$ is twice continuously differentiable for \ale $x\in\Omega$. 
			\end{enumerate}

		\item 
			$\varphi$ satisfies the boundedness and local Lipschitz conditions of order~$2$, \ie, there exists a constant $K > 0$ such that
			\begin{equation*}
				\abs[big]{D_y^\ell\varphi (x,0)}\le K \quad \text{for all } 0 \le \ell \le 2 \text{ and for \ale\ } x \in \Omega,
			\end{equation*} 
			and for every $M > 0$, there exists a Lipschitz constant $L(M) > 0$ such that
			\begin{equation*}
				\abs[big]{D_y^2\varphi(x,y_1) - D_y^2\varphi(x,y_2)}\le L(M) \, \abs{y_1-y_2}
			\end{equation*}
			holds for \ale $x\in \Omega$ and for all $\abs{y_i} \le M$, $i=1,2$.
	\end{enumeratearabic}
\end{assumption}

\cref{assumption:Nemytskii} implies the following properties for the Nemytskii operator $\Phi(y)(x) \coloneqq \varphi(x,y(x))$.

\begin{lemma}[\protect{\cite[Lemma~4.11, Lemma~4.12]{Troeltzsch:2010:1}}]\label{lemma:Nemytskii} \hfill 
	\begin{enumerate}
		\item 
			$\Phi$ is continuous in $L^\infty(\Omega)$. 
			Moreover, for all $r \in [1,\infty]$, we have 
			\begin{equation*}
				\norm{\Phi(y) - \Phi(z)}_{L^r(\Omega)} 
				\le 
				L(M) \, \norm{y-z}_{L^r(\Omega)}
			\end{equation*}
			for all $y,z\in L^\infty(\Omega)$ such that $\norm{y}_{L^\infty(\Omega)} \le M$ and $\norm{z}_{L^\infty(\Omega)} \le M.$

		\item  
			$\Phi$ is twice continuously Fréchet differentiable in $L^\infty(\Omega)$, and we have 
			\begin{equation*}
				\begin{aligned}
					\paren[auto](){\Phi'(y) \, h}(x) 
					&
					= 
					\varphi_y\paren[auto](){x,y(x)} h(x) 
					,
					\\
					\paren[auto](){\Phi''(y) \, [h_1,h_2]}(x) 
					&
					= 
					\varphi_{yy}\paren[auto](){x,y(x)} h_1(x) \, h_2(x) 
					.
				\end{aligned}
			\end{equation*}
			for \ale $x \in \Omega$ and $h, h_1, h_2 \in L^\infty(\Omega)$.
	\end{enumerate}
\end{lemma}

We now proceed to define the notion of weak solution of \eqref{eq:optimal_control_problem_Kirchhoff_equation}.
Since for any pair $(u,y)\in \Uad \times H^1(\Omega)$, $u + b \, \norm{\nabla y}_{L^2(\Omega)}^2$ is strictly positive, we can write the Kirchhoff equation \eqref{eq:optimal_control_problem_Kirchhoff_equation} in the form
\begin{equation}
	\label{eq:Kirchhoff_equation_divided_by_coefficient}
	-
	\laplace y 
	= 
	\frac{f}{u + b \, \norm{\nabla y}_{L^2(\Omega)}^2}
	.
\end{equation}
Here and in the following, we occasionally write $\norm{\cdot}$ instead of $\norm{\cdot}_{L^2(\Omega)}$. 
The $L^2(\Omega)$-inner product is denoted by $\inner{\cdot}{\cdot}$.
Moreover, we denote by $\cL(U,V)$ the space of bounded linear operators from $U$ to $V$.

Multiplication of \eqref{eq:Kirchhoff_equation_divided_by_coefficient} with a test function $v\in H^1_0(\Omega)$ and integration by parts yields the following definition. 
\begin{definition}\label{definition:weak_solution}
	A function $y \in H^1_0(\Omega)$ is called a weak solution of \eqref{eq:optimal_control_problem_Kirchhoff_equation} if it satisfies
	\begin{equation}\label{eq:weak_solution}
		\int_\Omega \nabla y \cdot \nabla v \dx 
		= 
		\int_\Omega \frac{f \, v}{u + b \, \norm{\nabla y}^2} \dx \quad \text{for all } v \in H^1_0(\Omega)  
		.
	\end{equation}
\end{definition}

The existence of a unique weak solution as well as its $W^{2,p}(\Omega)$-regularity, has been has shown in \cite[Theorem~2.2]{DelgadoFigueiredoGayteMoralesRodrigo:2017:1}.
Nevertheless, we briefly sketch the proof since its main idea is utilized again later on.
For a complete proof we refer the reader to \cite{DelgadoFigueiredoGayteMoralesRodrigo:2017:1}. 
\begin{theorem}
	\label{theorem:wellposedness_of_Kirchhoff_equation}
	For any $u \in \Uad$, there exists a unique weak solution $y \in H^1_0(\Omega)$ of the Kirchhoff problem \eqref{eq:optimal_control_problem_Kirchhoff_equation}. 
	Moreover, $y \in W^{2,p}(\Omega)$ holds for all $p \in [1,\infty)$, so it is also a strong solution.
\end{theorem}
\begin{proof}
	Suppose that $u \in \Uad$ and let $g \colon [0,\infty) \to \R$ be the function defined by 
	\begin{equation*}
		g(s) = s - \norm{\nabla y_s}^2
		,
	\end{equation*}
	where $y_s$ is the unique solution of the Poisson problem 
	\begin{equation*}
		\left\{
			\begin{aligned}
				-\laplace y_s & = \frac{f}{u + b \, s} & & \text{in } \Omega, \\
				y_s & = 0 & & \text{on } \partial\Omega.
			\end{aligned}
		\right.
	\end{equation*}
	A monotonicity argument can be used to show that $g$ has a unique root. 
	Since $y_s$ solves \eqref{eq:optimal_control_problem_Kirchhoff_equation} if and only if $g(s) = 0$ holds, the uniqueness of Kirchhoff equation is guaranteed. 
	Furthermore, due to the boundedness of $u$ from below, the right-hand side $f / (u + b \, s)$ of the Poisson problem above belongs to $L^\infty(\Omega)$. 
	Hence, by virtue of regularity results for the Poisson problem, $y \in W^{2,p}(\Omega)$ holds for any $p \in [1,\infty)$; see, \eg, \cite[Thm.~9.15]{GilbargTrudinger:1977:1}. 
\end{proof}

For the proof of existence of a globally optimal control of \eqref{eq:optimal_control_problem}, we show next that the control-to-state operator~$\controltostate \colon \Uad \to H^1_0(\Omega) \cap W^{2,p}(\Omega)$ is continuous. 
\begin{theorem}\label{theorem:control-to-state_map_continuity}
	The control-to-state map~$\controltostate$ is continuous from $\Uad$ (with the $L^2(\Omega)$-topology) into $H^1_0(\Omega) \cap W^{2,p}(\Omega)$ for all $p \in [1,\infty)$.
\end{theorem}
\begin{proof}
	The control-to-state map~$\controltostate \colon \Uad \to H^1_0(\Omega) \cap W^{2,p}(\Omega)$ is well-defined as a consequence of \cref{theorem:wellposedness_of_Kirchhoff_equation}.
	To show its continuity, let $\{u_n\} \subset \Uad$ be a sequence with $u_n \to u$ in $L^2(\Omega)$. 
	Set $y_n\coloneqq \controltostate(u_n)$, then we have the a-priori estimate 
	\makeatletter
	\begin{equation*}
		\begin{aligned}
			\norm{y_n}_{W^{2,p}(\Omega)} 
			&
			\le 
			c_1 \, \norm[auto]{\frac{f}{u_n + b \, \norm{\nabla y_n }^2}}_{L^p(\Omega)} 
			\le 
			c_2 \, \norm[auto]{\frac{f}{u_n + b \, \norm{\nabla y_n }^2}}_{L^\infty(\Omega)} 
			\ltx@ifclassloaded{svjour3}{%
				\\
				&
			}{}
			\le 
			c \, \norm[auto]{\frac{f}{u_a}}_{L^\infty(\Omega)}
			\le
			C
			.
		\end{aligned}
	\end{equation*}
	\makeatother
	From now on, suppose without loss of generality that $p \in [2,\infty)$ holds.
	Since $W^{2,p}(\Omega)$ is a reflexive Banach space and every bounded subset of a reflexive Banach space is weakly relatively compact, there exists a subsequence $y_n$, denoted by the same indices, satisfying $y_n \weakly \hat y$ in $W^{2,p}(\Omega)$. 
	The compactness of the embedding $W^{2,p}(\Omega) \embeds W^{1,p}(\Omega)$ implies the strong convergence $y_n \rightarrow \hat y$ in $W^{1,p}(\Omega)$ and thus $\norm{\nabla y_n} \to \norm{\nabla \hat y}$. 
	Moreover, $u_n \to u$ in $L^2(\Omega)$ implies the existence of a further subsequence $u_n$, still denoted by the same indices, with $u_n(x) \to u(x)$ for \ale $x \in \Omega$. 
	Consequently, 
	\begin{equation*}
		\frac{f}{u_n + b \, \norm{\nabla y_n }^2} \rightarrow \frac{f}{u + b \, \norm{\nabla\hat y}^2} 
		\quad 
		\text{\ale in } \Omega
		.
	\end{equation*}
	Since $\frac{f}{u_n + b \, \norm{\nabla y_n}^2}$ is dominated by $\frac{f}{u_a}$, we have 
	\begin{equation*}
		\abs[auto]{\frac{f}{u_n + b \, \norm{\nabla y_n }^2} 
			-
		\frac{f}{u + b \, \norm{\nabla\hat y}^2}}^p
		\le 
		\abs[auto]{\frac{2 f}{u_a}}^p
		.
	\end{equation*}
	By virtue of the dominated convergence theorem,
	\begin{equation*}
		- \laplace y_n 
		=
		\frac{f}{u_n + b \, \norm{\nabla y_n }^2} \rightarrow \frac{f}{u + b \, \norm{\nabla\hat y}^2}
		\quad 
		\text{in } L^p(\Omega)
		.
	\end{equation*}
	On the other hand, from $y_n \weakly \hat y$ in $W^{2,p}(\Omega)$, it follows that $\laplace y_n \weakly \laplace \hat y$ holds in $L^p(\Omega)$. 
	The uniqueness of the weak limit yields  
	\begin{equation*}
		- \laplace \hat y 
		=
		\frac{f}{u + b \, \norm{\nabla \hat y}^2}
	\end{equation*}
	and from the uniqueness of the solution of \eqref{eq:optimal_control_problem_Kirchhoff_equation} we obtain $\hat y = \controltostate(u)$.
	Therefore, $\laplace y_n \to \laplace \hat y$ holds in $L^p(\Omega)$ and thereby $y_n \to \hat y$ in $W^{2,p}(\Omega)$. 
	
	We note that we have proved that for any sequence $\{u_n\} \subset \Uad$ with $u_n \rightarrow u $ in  $L^2(\Omega)$ there exist a subsequence $\{u_{n}\}$, denoted by the same indices, so that $\controltostate(u_{n}) \rightarrow \controltostate(u)$ in $W^{2,p}(\Omega)$.
	Thus we can easily conclude  convergence of the entire sequence $\controltostate(u_n) \rightarrow \controltostate(u)$ in  $W^{2,p}(\Omega)$. 
	Indeed, if $\controltostate(u_n) \not \rightarrow \controltostate(u)$, then there exist $\delta > 0$ and a subsequences with indices~$n_k$ such that 
	\begin{equation*}
		\norm{\controltostate(u_{n_k}) - \controltostate(u)}_{W^{2,p}(\Omega)} 
		> 
		\delta
		\text{ for } k \to \infty
		.
	\end{equation*}
	Since $u_{n_k} \rightarrow u$ in  $L^2(\Omega)$, there exists a further subsequence $\{u_{n_{k_\ell}}\}$ such that $\controltostate(u_{n_{k_\ell}}) \rightarrow \controltostate(u)$, which is a contradiction.
	Consequently, we obtain $\controltostate(u_n) \rightarrow \controltostate(u)$ as claimed.
\end{proof}

The compact embedding $H^1(\Omega) \embeds L^2(\Omega)$ immediately leads to the following corollary.
\begin{corollary}\label{corollary:control-to-state_map_weak_strong_continuity}
	The control-to-state map~$\controltostate$ is weakly-strongly continuous from $\Uad \cap H^1(\Omega)$ (with the $H^1(\Omega)$-topology) into $H^1_0(\Omega) \cap W^{2,p}(\Omega)$ for all $p \in [1,\infty)$.
	That is, when $\{u_n\} \subset \Uad \cap H^1(\Omega)$ with $u_n \weakly u$ in $H^1(\Omega)$, then $\controltostate(u_n) \to \controltostate(u)$ in $W^{2,p}(\Omega)$. 
\end{corollary}

\begin{remark} \label{remark:necessity_of_the_upper_bound}
	An inspection of the proof of \cref{theorem:control-to-state_map_continuity} and \cref{corollary:control-to-state_map_weak_strong_continuity} shows that these results remain valid in the absence of an upper bound $u_b$ on the control. 
	However, the upper bound is used in the proof of existence of a globally optimal control in \cref{theorem:existence_of_global_minimizer} below.
\end{remark}

We can now address the existence of a global minimizer of \eqref{eq:optimal_control_problem}.
\begin{theorem}\label{theorem:existence_of_global_minimizer}
	Problem \eqref{eq:optimal_control_problem} possesses a globally optimal control $\bar u \in \Uad \cap H^1(\Omega)$ with associated optimal state $\bar y = \controltostate(\bar u) \in H^1_0(\Omega) \cap W^{2,p}(\Omega)$ for all $p \in [1,\infty)$.
\end{theorem}
\begin{proof}
	The proof follows the standard route of the direct method so we can be brief.
	\begin{enumerate}[label=Step~(\arabic*):,ref=Step~(\arabic*),leftmargin=*]
		\item 
			We show that the reduced cost functional 
			\begin{equation*}
				j(u)
				\coloneqq 
				\int_\Omega \Phi(\controltostate(u)) \dx 
				+ 
				\frac{\lambda}{2} \norm{u}_{H^1(\Omega)}^2
			\end{equation*}
			is bounded from below on the set $\Uad \cap H^1(\Omega)$.
			To this end, recall from the proof of \cref{theorem:control-to-state_map_continuity} that $\controltostate(\Uad)$ is bounded in $W^{2,p}(\Omega)$.
			Due to the embedding $W^{2,p}(\Omega) \embeds C(\overline \Omega)$ for $p > N / 2$, there exists $M > 0$ such that $\norm{\controltostate(u)}_{L^\infty(\Omega)} \le M$ holds for all $u \in \Uad$.
			From \cref{assumption:Nemytskii} we can obtain the estimate 
			\begin{align*}
				\abs{\varphi(x,\controltostate(u)(x))}
				& 
				\le 
				\abs{\varphi(x,0)} + \abs{\varphi(x,\controltostate(u)(x))-\varphi(x,0)} 
				\\
				& 
				\le 
				K + L(M) \, \abs{\controltostate(u)(x)}
				\le 
				K + L(M) \, M
				.
			\end{align*}
			This implies
			\begin{equation}\label{eq:existence_of_global_minimizer_boundedness_below}
				\int_\Omega \Phi(\controltostate(u)) \dx
				\ge 
				- \paren[big](){K + L(M) \, M} \, \abs{\Omega}
			\end{equation}
			for all $u \in \Uad$. 
			The assertion follows.

		\item
			We construct the tentative minimizer~$\bar u$.
			Since $j$ is bounded from below on $\Uad \cap H^1(\Omega)$, there exists a minimizing sequence $\{u_n\} \subset \Uad \cap H^1(\Omega)$ so that 
			\begin{equation*}
				j(u_n) \searrow \inf_{u \in \Uad} j(u) 
				\eqqcolon 
				\beta
				.
			\end{equation*}
			Consequently, there exists a subsequence, denoted by the same indices, such that $u_n \weakly \bar u$ in $H^1(\Omega)$.
			$\Uad$ is convex and closed in $H^1(\Omega)$ and therefore weakly closed in $H^1(\Omega)$, thus $\bar u \in \Uad \cap H^1(\Omega)$. 
			Now \cref{corollary:control-to-state_map_weak_strong_continuity} implies $\controltostate(u_n) \to \controltostate(\bar u)$ in $W^{2,p}(\Omega)$.

		\item
			It remains to show the global optimality of $\bar u$.
			Set $F(y) \coloneqq \int_\Omega \Phi(y) \dx$, thus $F$ is composed of a Nemytskii operator and a continuous linear integral operator from $L^1(\Omega)$ into $\R$. 
			By virtue of \cref{lemma:Nemytskii}, $\Phi$ is continuous in $L^\infty(\Omega)$. 
			Since $W^{2,p}(\Omega) \embeds L^\infty(\Omega)$ holds, $F \circ \controltostate$ is weakly-strongly continuous on $\Uad$ \wrt the topology of $L^2(\Omega)$.

			In summary, exploiting the weak sequential lower semicontinuity of $\norm{\cdot}_{H^1}$ we have
			\begin{align*}
				\beta 
				= 
				\lim_{n\to\infty} j(u_n)
				& 
				= 
				\lim_{n\to\infty} F(\controltostate(u_n)) + \frac{\lambda}{2} \liminf_{n\to\infty} \norm{u_n}_{H^1}^2 
				\\
				& 
				\ge
				F(\controltostate(\bar u)) + \frac{\lambda}{2} \norm{\bar u}_{H^1(\Omega)}^2
				= 
				j(\bar u)
				.
			\end{align*}
			By definition of $\beta$ and since $\bar u \in \Uad \cap H^1(\Omega)$, we therefore must have $\beta = j(\bar u)$.
			\ifthenelse{\isundefined{\qedhere}}{}{\qedhere}
	\end{enumerate}
\end{proof}

\section{Optimality System} 
\label{section:optimality_system}

In this section we address first-order necessary optimality conditions for local minimizers. 
With regard to an efficient numerical solution method in function spaces, we are aiming to arrive at an optimality system which is Newton differentiable.
Unfortunately, this is not the case for the first-order optimality system of our problem \eqref{eq:optimal_control_problem} due to the failure of the orthogonal projection \wrt $H^1(\Omega)$ onto $\Uad$ to be Newton differentiable.
We detail this issue in \cref{subsection:first-order_optimality}.
We therefore propose to relax and penalize the control constraints.
Notice that this is not straightforward since we need to ensure positivity of the relaxed control in the state equation.
We achieve the latter by a smooth cut-off function.
The optimality system of the penalized problem then turns out to be Newton differentiable, as we shall show in \cref{section:generalized_Newton_method}.

The material in this section is structured as follows.
In \cref{subsection:differentiability_control-to-state}, we prove the Fréchet differentiability of the control-to-state map.
We establish the system of first-order necessary optimality conditions for the original problem \eqref{eq:optimal_control_problem} in \cref{subsection:first-order_optimality}.
In \cref{subsection:MY_approximation} we introduce the penalty approximation and show that for any null sequence of penalty parameters, there exists a subsequence of global solutions to the corresponding penalized problems which converges weakly to a global solution of the original problem; see \cref{theorem:optimality_system_penalized_problem}.
\cref{subsection:first-order_optimality_penalized} addresses the system of first-order necessary optimality conditions for the penalized problem.

\subsection{Differentiability of the Control-to-State Map}
\label{subsection:differentiability_control-to-state}

In this subsection we show the Fréchet differentiability of the control-to-state map $\controltostate$ by means of the implicit function theorem. 
To verify the assumption of this theorem, we need the following result about the linearization of the Kirchhoff equation \eqref{eq:optimal_control_problem_Kirchhoff_equation}.
The proof idea is similar to \cref{theorem:wellposedness_of_Kirchhoff_equation} and the proof is omitted.
\begin{proposition}\label{proposition:linearized_Kirchhoff_equation}
	Suppose that $\hat u \in \Uad$ and $\hat y \in H^1_0(\Omega) \cap W^{2,p}(\Omega)$ is the associated unique solution of the Kirchhoff equation \eqref{eq:optimal_control_problem_Kirchhoff_equation} for any $p \in [1,\infty)$.
	Then, for any $g \in L^p(\Omega)$, the linearized problem
	\begin{equation}
		\label{eq:linearized_Kirchhoff_equation}
		\left\{
			\begin{aligned}
				-\paren[auto](){\hat u + b \, \norm{\nabla\hat y}^2} \laplace y
				-
				2 \, b \, \inner{\nabla \hat y}{\nabla y} \laplace \hat y
				&
				=
				g 
				& &
				\text{in } \Omega, 
				\\
				y
				&
				=
				0
				& &
				\text{on } \partial\Omega,
			\end{aligned}
		\right.
	\end{equation}
	has a unique solution~$y \in H^1_0(\Omega) \cap W^{2,p}(\Omega)$.
\end{proposition}

\begin{theorem}\label{theorem:Frechet_differentiability_control_to_state}
	The control-to-state operator~$\controltostate$
	\begin{equation*}
		\controltostate \colon \Uad \subset L^\infty(\Omega) \to H^1_0(\Omega) \cap W^{2,p}(\Omega)
	\end{equation*}
	is continuously Fréchet differentiable for all $p \in [1,\infty)$. 
\end{theorem}
\begin{proof}
	Suppose that $\hat u \in \Uad$ is arbitrary and that $\hat y \in H^1_0(\Omega) \cap W^{2,p}(\Omega)$ is the associated state.
	The map $E \colon \paren[auto](){H^1_0(\Omega) \cap W^{2,p}(\Omega)} \times L^\infty(\Omega) \to L^p(\Omega)$ defined by
	\begin{equation*}
		E(y,u)  
		\coloneqq
		- \paren[auto](){u + b \, \norm{\nabla y}^2} \laplace y - f
	\end{equation*}
	is continuously Fréchet differentiable with
	\begin{equation*}
		\begin{aligned}
			E'(\hat y , \hat u)(y , u) 
			=
			- \paren[big](){\hat u + b \, \norm{\nabla \hat y}^2} \laplace y
			- \paren[big](){u + 2 \, b \, \inner{\nabla \hat y}{\nabla y}} \laplace \hat y 
			.
		\end{aligned}
	\end{equation*}
	It remains to show that $E_y(\hat y,\hat u) \in \cL\paren[big](){H^1_0(\Omega) \cap W^{2,p}(\Omega), L^p(\Omega)}$ has a bounded inverse.
	To this end, consider
	\begin{equation}\label{eq:Frechet_differentiability_control_to_state_linearized_operator}
		E_y(\hat y,\hat u) \, y
		=
		- \paren[auto](){\hat u+b \, \norm{\nabla \hat y}^2} \laplace y 
		- 2 \, b \, \inner{\nabla \hat y}{\nabla y} \, \laplace \hat y
		.
	\end{equation} 
	The existence and uniqueness of $y \in H^1_0(\Omega) \cap W^{2,p}(\Omega)$ satisfying \eqref{eq:linearized_Kirchhoff_equation}, \ie, $E_y(\hat y,\hat u) \, y = g$, is established by virtue of \cref{proposition:linearized_Kirchhoff_equation}. 
	This implies the bijectivity of $E_y(\hat y,\hat u)$.
	The open mapping/continuous inverse theorem now yields that the inverse of $E_y(\hat y,\hat u)$ is continuous. 
	Notice that $E(y,u) = 0 \Leftrightarrow E(\controltostate(u),u) = 0$ holds for all $u \in \Uad$.
	Invoking the implicit function theorem, we obtain that $\controltostate$ is continuously differentiable in some $L^\infty(\Omega)$-neighborhood of $\hat u$.
	Since $\hat u \in \Uad$ was arbitrary, $\controltostate$ actually extends into an $L^\infty(\Omega)$-neighborhood of $\Uad$ and it is continuously differentiable there.
	Moreover, we obtain that $\delta y = \controltostate'(\hat u) \, \delta u$ satisfies $E_y(\hat y, \hat u) \, \delta y = - E_u(\hat y, \hat u) \, \delta u$, \ie,
	\begin{equation*}
		- \paren[big](){\hat u + b \, \norm{\nabla\hat y}^2} \laplace \delta y 
		- \paren[big](){\delta u + 2 \, b \, \inner{\nabla\hat y}{\nabla y}} \laplace \hat y
		=
		0
		.
		\ifthenelse{\isundefined{\qedhere}}{}{\qedhere}
	\end{equation*}
\end{proof}

\subsection{First-Order Optimality Conditions}
\label{subsection:first-order_optimality}

The optimality system can be derived by using the Lagrangian 
\begin{equation}
	\cL(y, u ,p) 
	\coloneq
	\int_\Omega \varphi(x, y) \dx 
	+
	\frac{\lambda}{2} \, \norm{u}_{H^1(\Omega)}^2 
	+
	\int_\Omega \nabla y \cdot \nabla p \dx 
	- 
	\int_\Omega \frac{f}{u + b \, \norm{\nabla y}^2} \, p \dx
\end{equation}
and taking the derivative with respect to the state and the control. 
In the first case, we obtain
\begin{equation*}
		\cL_y(y, u, p)\, \delta y 
		=
		\int_\Omega \varphi_y(x,y)\, \delta y \dx + \int_\Omega \nabla \delta y \cdot \nabla p \dx
		+
		\int_\Omega \frac{2\, b\, f\, p\, \inner{\nabla y}{\nabla \delta y}}{\paren[auto](){u + b \, \norm{\nabla y}^2}^2} \dx
\end{equation*}
for $\delta y \in H^1_0(\Omega) \cap W^{2,p}(\Omega)$.
Integration by parts yields
\makeatletter
\begin{equation*}
	\begin{aligned}
		\ltx@ifclassloaded{svjour3}{%
			\MoveEqLeft
			\cL_y(y, u, p)\, \delta y 
			\\
		}{%
			\cL_y(y, u, p)\, \delta y 
		}
		&
		=
		\int_\Omega \varphi_y(x,y)\, \delta y \dx
		+
		\int_\Omega \nabla \delta y \cdot \nabla p \dx
		+
		\inner[auto]{\nabla y\int_\Omega \frac{2\, b\, f\, p}{\paren[auto](){u + b \, \norm{\nabla y}^2}^2} \dx}{\nabla \delta y}
		\\
		&
		=
		\int_\Omega \varphi_y(x,y) \, \delta y \dx
		- 
		\int_\Omega \laplace p \, \delta y \dx
		-
		\inner[auto]{\laplace y\int_\Omega \frac{2\, b\, f\, p}{\paren[auto](){u + b \, \norm{\nabla y}^2}^2} \dx}{\delta y}
		.
	\end{aligned}
\end{equation*} 
\makeatother
Notice that $\cL_y(y, u, p)\, \delta y = 0$ for all $\delta y \in H^1_0(\Omega) \cap W^{2,p}(\Omega)$ represents the strong form of the adjoint equation, which reads
\begin{equation}
	\label{eq:adjoint_equation_strong_form}
	\left\{
		\begin{aligned}
			- \laplace p 
			- \laplace y \int_\Omega \frac{2 \, b \, f \, p}{\paren[auto](){u + b \, \norm{\nabla y}^2}^2}  \dx
			& 
			=
			- \varphi_y(x,y)
			& 
			& 
			\text{in } \Omega 
			,
			\\
			p 
			& 
			=
			0 
			& & 
			\text{on } \partial \Omega
			.
		\end{aligned}
	\right.
\end{equation}
We point out that \eqref{eq:adjoint_equation_strong_form} is again a nonlocal equation.
Given $u \in \Uad$ and $y \in H^1_0(\Omega) \cap W^{2,p}(\Omega)$, \eqref{eq:adjoint_equation_strong_form} has a unique solution $p \in H^1_0(\Omega) \cap W^{2,p}(\Omega)$.
This can be shown either by direct arguments as in \cref{theorem:wellposedness_of_Kirchhoff_equation}, or by exploiting that the bounded invertibility of $E_y$ implies that of its adjoint, see the proof of \cref{theorem:Frechet_differentiability_control_to_state}.

The derivative of the Lagrangian with respect to the control is given by
\begin{equation*}
	\begin{aligned}
		\cL_u(y , u, p) \, \delta u
		=
		\lambda \, \inner{u}{\delta u}_{H^1(\Omega)}
		+
		\int_\Omega \frac{f \, p}{\paren[auto](){u + b \, \norm{\nabla y}^2}^2} \, \delta u \dx
	\end{aligned}
\end{equation*}
for $\delta u \in L^\infty(\Omega) \cap H^1(\Omega)$.

It is now standard to derive the following system of necessary optimality conditions.
\begin{theorem}
	\label{theorem:optimality_conditions}
	Suppose that $(y,u) \in \paren[auto](){H^1_0(\Omega) \cap W^{2,p}(\Omega)} \times \paren[auto](){\Uad \cap H^1(\Omega)}$ is a locally optimal solution of problem \eqref{eq:optimal_control_problem} for any $p \in [1,\infty)$.
	Then there exists a unique adjoint state $p \in H^1_0(\Omega) \cap W^{2,p}(\Omega)$ for all $p \in [1,\infty)$ such that the following system holds:
	\begin{subequations}
		\label{eq:optimality_conditions}
		\begin{align}
			&
			\left\{
				\begin{aligned}
					- \laplace p 
					- \laplace y \int_\Omega \frac{2 \, b \, f \, p}{\paren[auto](){u + b \, \norm{\nabla y}^2}^2} \dx
					& 
					=
					- \varphi_y(x,y)
					& 
					& 
					\text{in } \Omega 
					,
					\\
					p 
					& 
					=
					0 
					& & 
					\text{on } \partial \Omega
					,
				\end{aligned}
			\right.
			\label{eq:optimality_conditions:adjoint}
			\\
			&
			\left\{
			\begin{aligned}
				& 
				\lambda \int_\Omega \nabla u \cdot \nabla(v - u) \dx
				+
				\int_\Omega \paren[auto](){\frac{f \, p}{\paren[auto](){u + b \, \norm{\nabla y}^2}^2} + \lambda \, u} \paren[auto](){v - u} \dx
				\ge
				0 
				\\
				& 
				\quad
				\text{for all } v \in \Uad \cap H^1(\Omega)
				,
			\end{aligned}
			\right.
			\label{eq:optimality_conditions:VI}
			\\
			&
			\left\{
			\begin{aligned}
				- \laplace y 
				&
				= 
				\frac{f}{u + b \, \norm{\nabla y}^2}
				& 
				& 
				\text{in } \Omega
				,
				\\
				y
				&
				=
				0
				&
				&
				\text{on } \partial \Omega
				.
			\end{aligned}
			\right.
			\label{eq:optimality_conditions:state}
		\end{align}
	\end{subequations}
\end{theorem}
Notice that \eqref{eq:optimality_conditions:VI} is a nonlinear obstacle problem for the control variable~$u$ originating from the bound constraints in $\Uad$ and the presence of the $H^1$-control cost term in the objective.
Unfortunately, this map is not known to be differentiable in the Newton sense.
In order to apply a generalized Newton method, we therefore relax and penalize the bound constraints via a quadratic penalty in the following section.
This is also known as Moreau-Yosida regularization of the indicator function pertaining to $\Uad$.

\subsection{Moreau-Yosida Penalty Approximation}
\label{subsection:MY_approximation}

The Moreau-Yosida penalty approximation of problem \eqref{eq:optimal_control_problem} consists of the following modifications.
\begin{enumeratearabic}
	\item 
		\label[modification]{item:remove_constraint}
		We remove the constraints $u_a \le u \le u_b$ from $\Uad$ and work with controls in $H^1(\Omega)$ which do not necessarily belong to $L^\infty(\Omega)$.

	\item
		\label[modification]{item:add_penalty}
		We add the penalty term $\frac{1}{2 \varepsilon} \int_\Omega \paren[auto](){u_a - u}_+^2 + \, \paren[auto](){u - u_b}_+^2 \dx$ to the objective.
		Here $v_+ = \max\{0,v\}$ is the positive part function and $\varepsilon > 0$ is the penalty parameter.

	\item
		\label[modification]{item:use_cutoff}
		We replace the control-to-state relation $y = S(u)$ by $y = \controltostate\paren[big](){u_a + \cutoff_\varepsilon(u - u_a)}$, where $\cutoff_\varepsilon$ is a family of monotone and convex $C^3$ approximations of the positive part function satisfying $\cutoff_\varepsilon(t) = t$ for $t > \varepsilon$, $\cutoff_\varepsilon(t) = 0$ for $t < - \varepsilon$ for some $\varepsilon > 0$ and $\cutoff_\varepsilon' \in [0,1]$ everywhere.
\end{enumeratearabic}
Notice that \cref{item:use_cutoff} is required since the control-to-state map~$S$ is guaranteed to be defined only for positive controls; compare \cref{theorem:control-to-state_map_continuity,remark:necessity_of_the_upper_bound}.
Therefore, we use $u_a + \cutoff_\varepsilon(u - u_a) \ge u_a$ as an effective control.
We now consider the following relaxed problem:
\begin{equation}
	\label{eq:optimal_control_problem_penalized}
  \tag{P$_\varepsilon$}
	\begin{aligned}
		\text{Minimize} \quad 
		&
		J_\varepsilon(y,u)
		\coloneqq
		J(y,u)
		+
		\frac{1}{2 \varepsilon} \int_\Omega \paren[auto](){u_a - u}_+^2 + \, \paren[auto](){u - u_b}_+^2 \dx
		\\
		\text{where } \quad
		&
		y = \controltostate\paren[auto](){u_a + \cutoff_\varepsilon(u - u_a)}
		\\
		\text{and} \quad
		&
		u \in H^1(\Omega) 
		.
	\end{aligned}
\end{equation}

The relation between \eqref{eq:optimal_control_problem_penalized} and the original problem \eqref{eq:optimal_control_problem} is clarified in the following theorem.
\begin{theorem}\label{theorem:optimality_system_penalized_problem} \hfill 
	\begin{enumerate}
		\item 
			\label[statement]{item:optimality_system_penalized_problem:1}
			For all $\varepsilon > 0$, problem \eqref{eq:optimal_control_problem_penalized} possesses a globally optimal solution $\paren(){\bar y_\varepsilon , \bar u_\varepsilon} \in \paren[auto](){H^1_0(\Omega) \cap W^{2,p}(\Omega)} \times H^1(\Omega)$ for all $p \in [1,\infty)$.

		\item 
			\label[statement]{item:optimality_system_penalized_problem:2}
			For any sequence $\varepsilon_n \searrow 0$, there is a subsequence of $\paren(){\bar y_{\varepsilon_n} , \bar u_{\varepsilon_n}}$ which converges weakly to some $\paren(){y^* , u^*}$ in $W^{2,p}(\Omega) \times H^1(\Omega)$.
			Moreover, $u^* \in \Uad$ holds and $\paren(){y^* , u^*}$ is a globally optimal solution of \eqref{eq:optimal_control_problem}.
	\end{enumerate}
\end{theorem} 
\begin{proof}
	\Cref{item:optimality_system_penalized_problem:1} can be proved in a straightforward manner using a similar procedure as in \cref{theorem:existence_of_global_minimizer}.
	The proof of \cref{item:optimality_system_penalized_problem:2} is divided into several steps.
	As in the proof of \cref{theorem:existence_of_global_minimizer}, we define $\beta$ to be the globally optimal value of the objective in \eqref{eq:optimal_control_problem}.
	Similarly, we let $\beta_\varepsilon$ denote the globally optimal value of the objective in \eqref{eq:optimal_control_problem_penalized}.
	Suppose that $\varepsilon_n \searrow 0$ is any sequence.

	\begin{enumerate}[label=Step~(\arabic*):,ref=Step~(\arabic*),leftmargin=*]
		\item
			\label[step]{item:optimality_system_penalized_problem:step1}
			We show that $\paren[auto]\{\}{\paren[auto](){\bar y_{\varepsilon_n} , \bar u_{\varepsilon_n}}}$ is bounded in $W^{2,p}(\Omega) \times H^1(\Omega)$. 

			Suppose that $(\bar y, \bar u)$ is a globally optimal solution of \eqref{eq:optimal_control_problem}.
			Owing to the definition of $\beta_\varepsilon$, we have
			\makeatletter
			\begin{equation}\label{eq:optimality_system_penalized_problem_1}
				\tag{$*$}
				\begin{aligned}
					\beta_\varepsilon 
					&
					\le
					J_\varepsilon(\bar y, \bar u)
					=
					J(\bar y , \bar u) 
					+ 
					\frac{1}{2 \varepsilon} \int_\Omega \paren[auto](){u_a - \bar u}_+^2 + \, \paren[auto](){\bar u - u_b}_+^2 \dx 
					\ltx@ifclassloaded{svjour3}{%
						\\
						&
					}{}
					= 
					J(\bar y , \bar u) 
					= 
					\beta
					.
				\end{aligned}
			\end{equation}
			\makeatother
			The next-to-last equality is true since $\bar u \in \Uad$ holds and therefore, the penalty term vanishes.
			Moreover, we obtain
			\makeatletter
			\begin{equation*}
				\begin{aligned}
					J(\bar y_{\varepsilon_n} , \bar u_{\varepsilon_n}) 
					&
					\le 
					J(\bar y_{\varepsilon_n} , \bar u_{\varepsilon_n}) 
					+
					\frac{1}{2 \varepsilon_n} \int_\Omega \paren[auto](){u_a - \bar u_{\varepsilon_n}}_+^2 + \, \paren[auto](){\bar u_{\varepsilon_n} - u_b}_+^2 \dx 
					\ltx@ifclassloaded{svjour3}{%
						\\
						&
					}{}
					=
					\beta_{\varepsilon_n}
					\le
					\beta
					,
				\end{aligned}
			\end{equation*}
			\makeatother
			where the last inequality follows from \eqref{eq:optimality_system_penalized_problem_1}.
			Since $\bar y_{\varepsilon_n} = \controltostate\paren[auto](){u_a + \cutoff(\bar u_{\varepsilon_n} - u_a)}$ holds, we obtain $\norm{\bar y_{\varepsilon_n}}_{W^{2,p}(\Omega)} \le C$ as in the proof of \cref{theorem:control-to-state_map_continuity}.
			Therefore, $\bar y_{\varepsilon_n}$ is also bounded in $C(\overline \Omega)$ and consequently, $\int_\Omega \varphi(x,\bar y_{\varepsilon_n}) \dx$ is bounded below, see \eqref{eq:existence_of_global_minimizer_boundedness_below}.
			Finally, 
			\begin{equation*}
				J(\bar y_{\varepsilon_n} , \bar u_{\varepsilon_n}) 
				= 
				\int_\Omega \varphi(x,\bar y_{\varepsilon_n}) \dx + \frac{\lambda}{2} \norm{\bar u_{\varepsilon_n}}_{H^1(\Omega)}^2
				\le 
				\beta
			\end{equation*}
			implies that $\norm{\bar u_{\varepsilon_n}}_{H^1(\Omega)}$ is bounded.

		\item
			\label[step]{item:optimality_system_penalized_problem:step2}
			From \cref{item:optimality_system_penalized_problem:step1} it follows that there exists a subsequence $\{\paren(){\bar y_{\varepsilon_n} , \bar u_{\varepsilon_n}}\}$, denoted with the same subscript, such that $\paren(){\bar y_{\varepsilon_n} , \bar u_{\varepsilon_n}} \weakly \paren(){y^* , u^*}$ in $W^{2,p}(\Omega) \times H^1(\Omega)$. 
			We show that $u^* \in \Uad$ holds.

			We have already shown that $\beta_{\varepsilon_n} \le \beta$ holds, therefore
			\begin{equation*}
				\int_\Omega \paren[auto](){u_a - \bar u_{\varepsilon_n}}_+^2 + \, \paren[auto](){\bar u_{\varepsilon_n} - u_b}_+^2 \dx 
				\le
				2 \varepsilon_n \paren[auto][]{\beta - J(\bar y_{\varepsilon_n} , \bar u_{\varepsilon_n})}
				.
			\end{equation*}
			Taking the $\limsup$ in this inequality as $n \to \infty$, we find 
			\makeatletter
			\begin{equation}\label{eq:optimality_system_penalized_problem_2}
				\tag{$**$}
				\begin{aligned}
					0
					&
					\le
					\limsup_{n \to \infty} \int_\Omega \paren[auto](){u_a - \bar u_{\varepsilon_n}}_+^2 + \, \paren[auto](){\bar u_{\varepsilon_n} - u_b}_+^2 \dx 
					\ltx@ifclassloaded{svjour3}{%
						\\
						&
					}{}
					\le 
					0 
					- 
					2 \, \liminf_{n\to\infty} \varepsilon_n \, J(\bar y_{\varepsilon_n}, \bar u_{\varepsilon_n})
					.
				\end{aligned}
			\end{equation}
			\makeatother
			From $\paren(){\bar y_{\varepsilon_n}, \bar u_{\varepsilon_n}} \weakly \paren(){y^* , u^*}$ in $W^{2,p}(\Omega) \times H^1(\Omega)$ we conclude $\bar u_{\varepsilon_n} \to u^*$ in $L^2(\Omega)$ and 
			\begin{equation}\label{eq:optimality_system_penalized_problem_3}
				\tag{$**$$*$}
				J(y^* , u^*) 
				\le
				\liminf_{n\to\infty} J(\bar y_{\varepsilon_n} , \bar u_{\varepsilon_n})
			\end{equation}
			as in the proof of \cref{theorem:existence_of_global_minimizer}.
			Passing with $n \to \infty$ in \eqref{eq:optimality_system_penalized_problem_2} yields 
			\begin{equation*}
				\int_\Omega \paren[auto](){u_a - u^*}_+^2 + \, \paren[auto](){u^* - u_b}_+^2 \dx 
				=
				0
			\end{equation*}
			and consequently, $u^* \in \Uad$ follows.

		\item
			\label{item:optimality_system_penalized_problem:step3}
			We show that $y^* = \controltostate(u^*)$ holds. 
			
			Since $\bar u_{\varepsilon_n} \to u^*$ in $L^2(\Omega)$, there exists a subsequence $\{\bar u_{\varepsilon_n}\}$, denoted with the same subscript, such that $\bar u_{\varepsilon_n}(x) - u_a(x) \to u^*(x) - u_a(x)$ for \ale $x \in \Omega$. Let $\delta > 0$ and $A_\delta$ be defined by  
			\begin{equation*}
				A_\delta
				= 
				\setDef{x \in \Omega}{u^*(x) - u_a(x) > \delta \text{ for } \ale\ x \in \Omega}
				.
			\end{equation*}
			We show that the admissibility of $u^*$ yields the following convergence:
			\begin{equation}
				\label{eq:optimality_system_penalized_problem_4}
				\cutoff_{\varepsilon_n}\paren[auto](){\bar u_{\varepsilon_n}(x) - u_a(x)} 
				\to
				u^*(x) - u_a(x)
				\quad \text{for } \ale\ x \in A_\delta
				.
			\end{equation}
			
			For $x \in A_\delta$, we have $\bar u_{\varepsilon_n}(x) - u_a(x) > \delta / 2$ for $n$ sufficiently large. 
			Therefore for $\varepsilon_n  < \delta / 2$ we can conclude 
			\begin{equation*}
				\cutoff_{\varepsilon_n}\paren[auto](){\bar u_{\varepsilon_n}(x) - u_a(x)} 
				=
				\bar u_{\varepsilon_n}(x) - u_a(x)
				.
			\end{equation*}
			We note that if $u^*(x) = u_a(x)$, we find $\bar u_{\varepsilon_n}(x) - u_a(x) \to u^*(x) - u_a(x) = 0$.   
			From 
			\begin{equation}
				\label{eq:optimality_system_penalized_problem_5}
				\abs[auto]{\cutoff_{\varepsilon_n}\paren[auto](){\bar u_{\varepsilon_n}(x) - u_a(x)}}
				\le
				\varepsilon_n + \abs[auto]{\bar u_{\varepsilon_n}(x) - u_a(x)}
				\quad \text{for } \ale\ x \in \Omega
			\end{equation}
			it follows that $\cutoff_{\varepsilon_n}\paren[auto](){\bar u_{\varepsilon_n}(x) - u_a(x)} \to 0$. 
			Since $\delta$ is arbitrary, \eqref{eq:optimality_system_penalized_problem_4} is valid for almost all $x$ with $u^*(x) \ge u_a(x)$.
			
			From \eqref{eq:optimality_system_penalized_problem_5} and \eqref{eq:optimality_system_penalized_problem_4} and by the dominated convergence theorem we obtain
			\begin{equation*}
				\cutoff_{\varepsilon_n}\paren[auto](){\bar u_{\varepsilon_n} - u_a} 
				\to
				u^* - u_a \quad \text{in } L^2(\Omega)
				.
			\end{equation*}
			The continuity of $\controltostate$ on $\Uad$ \wrt the $L^2(\Omega)$-topology now implies
			\begin{equation*}
				\bar y_{\varepsilon_n} 
				=
				\controltostate \paren[auto](){u_a + \cutoff_{\varepsilon_n}\paren[auto](){\bar u_{\varepsilon_n}(x) - u_a(x)}} 
				\to
				\controltostate \paren[auto](){u^*(x)}
				.
			\end{equation*}
			From \cref{item:optimality_system_penalized_problem:step2} we have the weak convergence of $\bar y_{\varepsilon_n}$ to $y^*$. 
			The uniqueness of the weak limit shows $y^* = \controltostate(u^*)$.
			
		\item
			Since $J(\bar y_{\varepsilon_n} , \bar u_{\varepsilon_n}) \le \beta$ holds, we obtain $J(y^* , u^*) \le \beta$ by invoking \eqref{eq:optimality_system_penalized_problem_3}. 
			Moreover, since $(y^*,u^*)$ is admissible for \eqref{eq:optimal_control_problem}, the definition of $\beta$ implies $J(y^*,u^*) = \beta$, which completes the proof.
			\ifthenelse{\isundefined{\qedhere}}{}{\qedhere}
	\end{enumerate}
\end{proof}

\subsection{First-Order Optimality Conditions for the Penalized Problem}
\label{subsection:first-order_optimality_penalized}

The derivation of optimality conditions for \eqref{eq:optimal_control_problem_penalized} proceeds along the same lines as in \cref{subsection:first-order_optimality} and the details are omitted. 
For simplicity, we drop the index $\cdot_\varepsilon$ from now on and denote states, controls, and associated adjoint states by $(y,u,p)$.
We obtain the following regularized system of necessary optimality conditions.
\begin{theorem}
	\label{theorem:optimality_conditions_regularized}
	Suppose that $(y,u) \in \paren[auto](){H^1_0(\Omega) \cap W^{2,p}(\Omega)} \times H^1(\Omega) $ is a locally optimal solution of problem \eqref{eq:optimal_control_problem_penalized} for any $p \in [1,\infty)$.
	Then there exists a unique adjoint state $p \in H^1_0(\Omega) \cap W^{2,p}(\Omega)$ for all $p \in [1,\infty)$ such that the following system holds:
	\begin{subequations}
		\label{eq:optimality_conditions_regularized}
		\begin{align}
			&
			\left\{
				\begin{aligned}
					- \laplace p 
					- \laplace y \int_\Omega \frac{2 \, b \, f \, p}{\paren[auto](){u_a + \cutoff_\varepsilon\paren[auto](){u - u_a} + b \, \norm{\nabla y}^2}^2} \dx
					& 
					=
					- \varphi_y(x,y)
					& 
					& 
					\text{in } \Omega 
					,
					\\
					p 
					& 
					=
					0 
					& & 
					\text{on } \partial \Omega
					,
				\end{aligned}
			\right.
			\label{eq:optimality_conditions_regularized:adjoint}
			\\
			&
			\left\{
			\begin{aligned}
				& 
				\lambda \int_\Omega \nabla u \cdot \nabla v \dx
				+
				\int_\Omega \paren[auto](){\frac{f \, p \, \cutoff_\varepsilon'(u - u_a)}{\paren[auto](){u_a + \cutoff_\varepsilon\paren[auto](){u - u_a} + b \, \norm{\nabla y}^2}^2} + \lambda \, u} v \dx
				\\
				& 
				\quad
				- 
				\frac{1}{\varepsilon} \int_\Omega \paren[big](){\paren[auto](){u_a - u}_+ - \paren[auto](){u - u_b}_+} \, v \dx
				=
				0 
				\quad
				\text{for all } v \in H^1(\Omega) 
				,
			\end{aligned}
			\right.
			\label{eq:optimality_conditions_regularized:VI}
			\\
			&
			\left\{
			\begin{aligned}
				- \laplace y 
				&
				= 
				\frac{f}{u_a + \cutoff_\varepsilon\paren[auto](){u - u_a} + b \, \norm{\nabla y}^2}
				& 
				& 
				\text{in } \Omega
				,
				\\
				y
				&
				=
				0
				&
				&
				\text{on } \partial \Omega
				.
			\end{aligned}
			\right.
			\label{eq:optimality_conditions_regularized:state}
		\end{align}
	\end{subequations}
\end{theorem}

\begin{corollary}
	\label{corollary:regularity_of_optimal_controls}
	The terms 
	\begin{equation*}
		\paren[auto](){\frac{f \, p \, \cutoff_\varepsilon'(u - u_a)}{\paren[auto](){u_a + \cutoff_\varepsilon\paren[auto](){u - u_a} + b \, \norm{\nabla y}^2}^2} + \lambda \, u} 
		- 
		\frac{1}{\varepsilon} \paren[big](){\paren[auto](){u_a - u}_+ - \paren[auto](){u - u_b}_+}
	\end{equation*}
	in \eqref{eq:optimality_conditions_regularized:VI} belong to $L^\infty(\Omega)$ and therefore, any locally optimal control of \eqref{eq:optimal_control_problem_penalized} belongs to $W^{2,p}(\Omega)$ for any $p \in [1,\infty)$.
\end{corollary}
\begin{proof}
	We only elaborate on the case $N = 3$ since the cases $N \in \{1,2\}$ are similar.
	We first consider the numerator of the first term.
	Here $f \in L^\infty(\Omega)$ holds by \cref{assumption:general} and $p \in L^\infty(\Omega)$ by virtue of the embedding $W^{2,p}(\Omega) \embeds L^\infty(\Omega)$ for $p > 3/2$.
	Moreover, $\cutoff_\varepsilon'$ maps into $[0,1]$ and therefore $\cutoff_\varepsilon'(u - u_a)$ belongs to $L^\infty(\Omega)$ as well.
	The denominator is bounded below by $u_a$, and therefore, the first term belongs to $L^\infty(\Omega)$.
	\\
	The second term, $\frac{1}{\varepsilon} \paren[big](){\paren[auto](){u_a - u}_+ - \paren[auto](){u - u_b}_+}$, belongs to $L^6(\Omega)$ due to the embedding $H^1(\Omega) \embeds L^6(\Omega)$.
	Inserting this into \eqref{eq:optimality_conditions_regularized:VI} with the differential operator $\lambda \, (-\laplace + \id)$ and the remaining terms on the right-hand side shows $u \in W^{2,6}(\Omega)$, which in turn embeds into $L^\infty(\Omega)$.
	\\
	Repeating this procedure one more time implies $u \in W^{2,p}(\Omega)$.
\end{proof}

\section{Generalized Newton Method}
\label{section:generalized_Newton_method}

In this section we show that the optimality system \eqref{eq:optimality_conditions_regularized} of the penalized problem is differentiable in a generalized sense, referred to as Newton differentiability.
This allows us to formulate a generalized Newton method.
Due to its similarity with the concept of semismoothness, see \cite{Ulbrich:2011:1}, such methods are sometimes referred to as a semismooth Newton method.

\begin{definition}[\protect{\cite[Definition~1]{HintermuellerItoKunisch:2002:1}, \cite[Definition~8.10]{ItoKunisch:2008:1}}]
	\label{definition:Newton_differentiability}
	Let $D$ be an open subset of a Banach space~$X$. 
	The mapping $F \colon D \subset X \to Y$ is called Newton differentiable on the open subset $V \subset D$ if there exists a map $G \colon V \to \cL(X,Y)$ such that, for every $x \in V$,
	\begin{equation*}
		\lim_{h \to 0} \frac{1}{\norm{h}_X} \norm[auto]{F(x + h)- F(x) - G(x + h) h}_Y 
		=
		0
		.
	\end{equation*} 
	In this case $G$ is said to be a Newton derivative of $F$ on $V$.
\end{definition}

We formulate the optimality system \eqref{eq:optimality_conditions_regularized} in terms of an operator equation $F = 0$ where
\begin{equation}
	\label{eq:optimality_conditions_regularized_operator}
	F \colon X \coloneqq \paren[auto](){W^{2,p}(\Omega) \cap H_0^1(\Omega)} \times W^{2,p}_\diamond(\Omega) \times \paren[auto](){W^{2,p}(\Omega) \cap H_0^1(\Omega)} 
	\to
	L^p(\Omega)^3 \eqqcolon Y
\end{equation}
and $p \in [\max\{1,N/2\},\infty)$ is arbitrary but fixed. Here $W^{2,p}_\diamond(\Omega)$ is defined as
\begin{equation*}
    W^{2,p}_\diamond(\Omega)
	\coloneqq
	\setDef[auto]{u \in W^{2,p}(\Omega)}{\frac{\partial u}{\partial n} = 0 \text{ on } \partial\Omega}
	.
\end{equation*}
The component $F_1$ represents the adjoint equation \eqref{eq:optimality_conditions_regularized:adjoint} in strong form, \ie,
\begin{equation*}
	F_1(y,u,p)
	=
	- \laplace p 
	- \laplace y \int_\Omega \frac{2 \, b \, f \, p}{\paren[auto](){u_a + \cutoff_\varepsilon\paren[auto](){u - u_a} + b \, \norm{\nabla y}^2}^2} \dx
	+ \varphi_y(x,y)
	.
\end{equation*}
The continuous Fréchet differentiability of $F_1$ is a standard result, which uses \cref{lemma:Nemytskii} and the embedding $W^{2,p}(\Omega) \embeds L^\infty(\Omega)$.
The directional derivative is given by 
\makeatletter
\begin{equation*}
	\begin{aligned}
		\MoveEqLeft
		F_1'(y,u,p) \, (\delta y, \delta u, \delta p)
		\\
		&
		=
		- \laplace \delta p 
		- \laplace \delta y \int_\Omega \frac{2 \, b \, f \, p}{\paren[auto](){u_a + \cutoff_\varepsilon\paren[auto](){u - u_a} + b \, \norm{\nabla y}^2}^2} \dx
		\ltx@ifclassloaded{svjour3}{%
			\\
			&
			\quad
		}{}
		- \laplace y \int_\Omega \frac{2 \, b \, f \, \delta p}{\paren[auto](){u_a + \cutoff_\varepsilon\paren[auto](){u - u_a} + b \, \norm{\nabla y}^2}^2} \dx 
		\\
		&
		\quad
		+ \laplace y \int_\Omega \frac{4 \, b \, f \, p \paren[auto](){\cutoff_\varepsilon'(u - u_a) \, \delta u + 2 \, b \, \inner{\nabla y}{\nabla \delta y}}}{\paren[auto](){u_a + \cutoff_\varepsilon\paren[auto](){u - u_a} + b \, \norm{\nabla y}^2}^3} \dx
		+ \varphi_{yy}(x,y) \, \delta y
		.
	\end{aligned}
\end{equation*}
\makeatother

Similarly, $F_3$ represents the state equation \eqref{eq:optimality_conditions_regularized:state}, \ie,
\begin{equation*}
	F_3(y,u,p)
	=
	- \laplace y 
	- \frac{f}{u_a + \cutoff_\varepsilon\paren[auto](){u - u_a} + b \, \norm{\nabla y}^2}
\end{equation*}
and its continuous Fréchet derivative is given by
\begin{equation*}
	F_3'(y,u,p) (\delta y, \delta u, \delta p)
	=
	- \laplace \delta y 
	+ \frac{f \paren[auto][]{\cutoff_\varepsilon'(u - u_a) \, \delta u + 2 \, b \, \inner{\nabla y}{\nabla \delta y}}}{\paren[auto](){u_a + \cutoff_\varepsilon\paren[auto](){u - u_a} + b \, \norm{\nabla y}^2}^2}
	.
\end{equation*}
Finally, in order to define $F_2$ we integrate \eqref{eq:optimality_conditions_regularized:VI} by parts, which is feasible due to \cref{corollary:regularity_of_optimal_controls}.
This results in the equivalent formulation $F_2 = 0$, where
\makeatletter
\begin{equation*}
	\begin{aligned}
		F_2(y,u,p)
		&
		=
		- \lambda \laplace u 
		+
		\frac{f \, p \, \cutoff_\varepsilon'(u - u_a)}{\paren[auto](){u_a + \cutoff_\varepsilon\paren[auto](){u - u_a} + b \, \norm{\nabla y}^2}^2} 
		+ \lambda \, u	
		\ltx@ifclassloaded{svjour3}{%
			\\
			&
			\quad
		}{}
		- 
		\frac{1}{\varepsilon} \paren[auto](){\max\{u_a - u, 0\} - \max\{u - u_b, 0\}}
		, 
	\end{aligned}
\end{equation*}
\makeatother
and the boundary conditions $\frac{\partial u}{\partial n} = 0$, which are included in the definition of $W^{2,p}_\diamond(\Omega)$.

In order to establish the Newton differentiability of $F_2$, we invoke the following classical result.
\begin{theorem}[\protect{\cite[Propoposition~4.1]{HintermuellerItoKunisch:2002:1}, \cite[Example~8.14]{ItoKunisch:2008:1}}]
	\label{theorem:Newton_differentiability_max_function}
	The mapping $\max\{0,\cdot\} \colon L^q(\Omega) \to L^p(\Omega)$, $1 \le p < q \le \infty$ is Newton differentiable on $L^q(\Omega)$ with generalized derivative 
	\begin{equation*}
		G_{\max}
		\colon
		L^q(\Omega) \to \cL(L^q(\Omega) , L^p(\Omega))
	\end{equation*}
	given by 
	\begin{equation*}
		G_{\max}(u) \, \delta u
		=
		\begin{cases}
			\delta u(x), & \text{where } u(x) > 0
			,
			\\
			0, & \text{where } u(x) \le 0
			.
		\end{cases}
	\end{equation*}
\end{theorem}

Using \cref{theorem:Newton_differentiability_max_function} and the embedding $W^{2,p} \embeds L^\infty(\Omega)$, it follows that $F_2$ is Newton differentiable on the entire space~$X$ with generalized derivative
\begin{equation*}
	\begin{aligned}
		\MoveEqLeft
		G_2(y,u,p)(\delta y, \delta u, \delta p)
		\\
		&
		=
		- \lambda \laplace \delta u 
		+ \frac{f \, \delta p \, \cutoff_\varepsilon'(u - u_a) + f \, p \, \cutoff_\varepsilon''(u - u_a) \, \delta u}{\paren[auto](){u_a + \cutoff_\varepsilon\paren[auto](){u - u_a} + b \, \norm{\nabla y}^2}^2} 
		\\
		&
		\quad
		- \frac{2 \, f \, p \, \cutoff_\varepsilon'(u - u_a) \paren[auto][]{\cutoff_\varepsilon'(u - u_a) \, \delta u + 2 \, b \, \inner{\nabla y}{\nabla \delta y}}}{\paren[auto](){u_a + \cutoff_\varepsilon\paren[auto](){u - u_a} + b \, \norm{\nabla y}^2}^3} 
		+ \lambda \, \delta u 
		+ \frac{1}{\varepsilon} \chi_{A(u)} \delta u
		.
	\end{aligned}
\end{equation*}
Here $\chi_A$ stands for the indicator function of the set 
\begin{equation*}
	A(u)
	=
	\setDef{x \in \Omega}{u_a - u \ge 0  \text{ or } u - u_b \ge 0}
	.
\end{equation*}

We are now in a position to state a basic generalized Newton method; see \cref{algorithm:semismooth_Newton_method}.
Following well-known arguments, we can show its local well-posedness and superlinear convergence to local minimizers satisfying second-order sufficient conditions.
We refrain from repeating the details and refer the interested reader to, \eg, \cite[Ch.~7]{ItoKunisch:2008:1}, \cite[Ch.~2.4--2.5]{HinzePinnauUlbrichUlbrich:2009:1} and \cite[Ch.~10]{Ulbrich:2011:1}.
It is also possible to globalize the method using a line seach approach; see, \eg, \cite{HinzeVierling:2012:1}.

\begin{algorithm}[Basic semismooth Newton method for the solution of problem \eqref{eq:optimal_control_problem_penalized}] \hfill
	\label{algorithm:semismooth_Newton_method}
	\begin{algorithmic}[1]
		\Require initial guess $(y_0,u_0,p_0) \in X$ 
		\Ensure approximate stationary point of \eqref{eq:optimal_control_problem_penalized}
		\State Set $k \coloneqq 0$
		\While{not converged}
		\State Determine the active set $A(u_k)$
		\State Solve the Newton system
		\begin{equation}
			\label{eq:Newton_system}
			\begin{aligned}
				G_1(y_k,u_k,p_k) (\delta y, \delta u, \delta p)
				&
				=
				- F_1(y_k,u_k,p_k)
				\\
				G_2(y_k,u_k,p_k) (\delta y, \delta u, \delta p)
				&
				=
				- F_2(y_k,u_k,p_k)
				\\
				G_3(y_k,u_k,p_k) (\delta y, \delta u, \delta p)
				&
				=
				- F_3(y_k,u_k,p_k)
				\\
			\end{aligned}
		\end{equation}
		\State Update the iterates by setting
		\begin{equation*}
			y_{k+1} 
			\coloneqq
			y_k + \delta y
			,
			\quad
			u_{k+1} 
			\coloneqq
			u_k + \delta u
			,
			\quad
			p_{k+1} 
			\coloneqq
			p_k + \delta p
		\end{equation*}
		\State Set $k \coloneqq k+1$
		\EndWhile
	\end{algorithmic}
\end{algorithm}

An appropriate criterion for the convergence of \cref{algorithm:semismooth_Newton_method} is the smallness of $\norm{F_1(y_k,u_k,p_k)}_{L^p(\Omega)}$, $\norm{F_2(y_k,u_k,p_k)}_{L^p(\Omega)}$ and $\norm{F_3(y_k,u_k,p_k)}_{L^p(\Omega)}$, either in absolute terms or relative to the initial values.

\begin{remark}
	\label{remark:convex_domains}
	We remark that all previous results can be generalized to convex domains $\Omega \subset \R^N$ where $1 \le N \le 3$.
	In this case, we can invoke the $H^2$-regularity result for the Poisson problem on convex domains from \cite[Thm.~3.2.1.3]{Grisvard:1985:1} in the proof of \cref{theorem:wellposedness_of_Kirchhoff_equation}.
	Consequently, we have to replace $p \in [1,\infty)$ by $p = 2$ in \cref{theorem:wellposedness_of_Kirchhoff_equation} and all subsequent results.
	The requirement $N \le 3$ ensures the validity of the embedding $H^2(\Omega) \embeds C(\overline \Omega)$.
\end{remark}

\section{Discretization and Implementation}
\label{section:discretization}

In this section we address the discretization of the relaxed optimal control problem \eqref{eq:optimal_control_problem_penalized}.
We then follow a discretize--then optimize approach and derive the associated discrete optimality system, as well as a discrete version of the generalized Newton method.
In order to simplify the implementation, we employ the original control-to-state map $y = \controltostate(u)$.
In other words, we choose $\cutoff_\varepsilon = \id$ in \eqref{eq:optimal_control_problem_penalized}, which no longer approximates the positive part function.
Consequently, the controls appearing in the control-to-state map are no longer guaranteed to be bounded below by $u_a$.
This simplification is justified as long as the control iterates still permit the state equation to be uniquely solvable, or rather its linearized counterpart appearing in the generalized Newton method.

Our discretization method of choice is the finite element method.
We employ piecewise linear, globally continuous finite elements on geometrically conforming triangulations of the domain~$\Omega$.
More precisely, we use the space
\begin{equation*}
	V_h 
	\coloneqq
	\setDef{v \in H^1(\Omega) \cap C(\overline \Omega)}{v \text{ is linear on all triangles}}
	\subset H^1(\Omega)
\end{equation*}
to discretize the control, the state and adjoint state variables.
We use the usual Lagrangian basis and refer to the basis functions as $\{\varphi_j\}$, where $j = 1, \ldots, N_V$ and $N_V$ denotes the number of vertices in the mesh.
The coefficient vector, \eg, for the discrete control variable~$u \in V_h$, will be denoted by $\tu$, so we have
\begin{equation*}
	u 
	= 
	\sum_{j=1}^{N_V} \tu_j \varphi_j
	.
\end{equation*}
In order to formulate the discrete optimal control problem, we introduce the mass and stiffness matrices $\tM$ and $\tK$ as follows:
\begin{equation*}
	\tM_{ij} 
	=
	\int_\Omega \varphi_i \, \varphi_j \dx
	\quad \text{and} \quad
	\tK_{ij} 
	=
	\int_\Omega \nabla \varphi_i \cdot \nabla \varphi_j \dx
	.
\end{equation*}
We also make use of the diagonally lumped mass matrix $\tMlumped$ with entries $\tMlumped_{ii} = \sum_{j=1}^{N_V} \tM_{ij}$.
Suppose that the right-hand side~$f$ and coefficient~$b$ have been discretized and represented by their coefficient vectors~$\tf$ and $\tb$ in $V_h$.
Using the lumped mass matrix, the weak formulation \eqref{eq:weak_solution} of the state equation can be written in preliminary discrete form as 
\begin{equation*}
	\tK \, \ty
	=
	\tMlumped \paren[auto][]{\frac{\tf_i}{\tu_i + \tb_i (\ty^\transp \tK \, \ty)}}_{i=1}^{N_V}
	.
\end{equation*}
In order to incorporate the Dirichlet boundary conditions, we introduce the boundary projector $\tP_\Gamma$.
This is a diagonal $N_V \times N_V$-matrix which has ones along the diagonal in entries pertaining to boundary vertices, and zeros otherwise.
We also introduce the interior projector $\tP_\Omega \coloneqq \id - \tP_\Gamma$.
We can thus state the discrete form of the state equation \eqref{eq:weak_solution} as
\begin{equation}
	\label{eq:state_equation_discrete_with_vectors}
	\tP_\Omega \tK \, \ty
	-
	\tP_\Omega \tMlumped \paren[auto][]{\frac{\tf_i}{\tu_i + \tb_i (\ty^\transp \tK \, \ty)}}_{i=1}^{N_V}
	+ 
	\tP_\Gamma \ty
	=
	\bnull
\end{equation}
In order to simplify the notation, we introduce further diagonal matrices
\begin{equation*}
	\tF 
	\coloneqq
	\diag(\tf)
	,
	\quad
	\tB 
	\coloneqq
	\diag(\tb)
	\quad \text{and} \quad
	\tD(\ty,\tu)
	\coloneqq
	\diag(\tu) + (\by^\transp \tK \, \ty) \, \tB
	.
\end{equation*}
Using these matrices, we can write \eqref{eq:state_equation_discrete_with_vectors} more compactly as
\begin{equation}
	\label{eq:state_equation_discrete}
	e(\ty,\tu)
	\coloneqq
	\tP_\Omega \tK \, \ty
	-
	\tP_\Omega \tMlumped \, \tF \, \tD(\ty,\tu)^{-1} \bone
	+ 
	\tP_\Gamma \ty
	=
	\bnull
	,
\end{equation}
where $\bone$ and $\bnull$ denote column vectors of all ones and all zeros, respectively.

To be specific, we focus on a tracking-type objective and choose $\varphi(x,y) = \frac{1}{2} (y-y_d)^2$ in \eqref{eq:optimal_control_problem_penalized}.
Consequently, the objective is discretized as
\begin{multline}
	J(\ty,\tu)
	=
	\frac{1}{2} (\ty - \ty_d)^\transp \tM (\ty - \ty_d)
	+
	\frac{\lambda}{2} \tu^\transp (\tK + \tM) \tu
	\\
	+
	\frac{1}{2 \varepsilon} (\tu_a - \tu)_+^\transp \tMlumped (\tu_a - \tu)_+
	+
	\frac{1}{2 \varepsilon} (\tu - \tu_b)_+^\transp \tMlumped (\tu - \tu_b)_+
	\label{eq:objective_discrete}
\end{multline}
and the Lagrangian of our discretized problem becomes
\begin{multline}
	\cL(\ty,\tu,\tp)
	=
	\frac{1}{2} (\ty - \ty_d)^\transp \tM (\ty - \ty_d)
	+
	\frac{\lambda}{2} \tu^\transp (\tK + \tM) \tu
	\\
	+
	\frac{1}{2 \varepsilon} (\tu_a - \tu)_+^\transp \tMlumped (\tu_a - \tu)_+
	+
	\frac{1}{2 \varepsilon} (\tu - \tu_b)_+^\transp \tMlumped (\tu - \tu_b)_+
	\\
	+ 
	\tp^\transp \tP_\Omega \tK \, \ty
	-
	\tp^\transp \tP_\Omega \tMlumped \, \tF \, \tD(\ty,\tu)^{-1} \bone
	+ 
	\tp^\transp \tP_\Gamma \ty
	.
	\label{eq:Lagrangian_discrete}
\end{multline}
Before we state the first- and second-order derivatives of the Lagrangian, we address the nonlinear term $\tD(\ty,\tu)^{-1}$ first.
We obtain
\makeatletter
\begin{equation*}
	\begin{alignedat}{2}
		\frac{\d}{\d \ty} \tD(\ty,\tu)^{-1} \delta \ty
		&
		=
		- 2 \, (\ty^\transp \tK \, \delta \ty) \, \tB \, \tD(\ty,\tu)^{-2} 
		&
		&
		\ltx@ifclassloaded{svjour3}{%
			\\
			\text{ and thus }
		}{%
			\text{ and thus }
		}
		\frac{\d}{\d \ty} \tD(\ty,\tu)^{-1} \bone
		\ltx@ifclassloaded{svjour3}{&}{}
		=
		-2 \, \tB \, \tD(\ty,\tu)^{-2} \bone \, \ty^\transp \tK
		,
		\\
		\frac{\d}{\d \tu} \tD(\ty,\tu)^{-1} \delta \tu
		&
		=
		- \tD(\ty,\tu)^{-2} \diag(\delta \tu)
		&
		&
		\ltx@ifclassloaded{svjour3}{%
			\\
			\text{ and thus }
		}{%
			\text{ and thus }
		}
		\frac{\d}{\d \tu} \tD(\ty,\tu)^{-1} \bone
		\ltx@ifclassloaded{svjour3}{&}{}
		=
		- \tD(\ty,\tu)^{-2}
		.
	\end{alignedat}
\end{equation*}
\makeatother
Therefore, the first-order derivatives of $\cL$ (written as column vectors) are given by
\makeatletter
\begin{subequations}
	\label{eq:first-order_derivatives_Lagrangian_discrete}
	\begin{align}
		\cL_\by(\by,\bu,\bp)
		&
		=
		\tM (\ty - \ty_d)
		+
		\tK \, \tP_\Omega \tp 
		\ltx@ifclassloaded{svjour3}{%
			\notag
			\\
			&
			\quad
		}{}
		+
		2 \, \tK \, \ty \, \bone^\transp \tD(\ty,\tu)^{-2} \tB \, \tF \, \tMlumped \, \tP_\Omega \tp
		+ 
		\tP_\Gamma \tp 
		,
		\label{eq:first-order_derivative_Lagrangian_wrt_y_discrete}
	\end{align}
	and
	\begin{multline}
		\cL_\bu(\by,\bu,\bp)
		=
		\lambda \, (\tK + \tM) \, \tu
		-
		\frac{1}{\varepsilon} \tD_{A_-}(\tu) \, \tMlumped \, \tD_{A_-}(\tu) \, (\tu_a - \tu)
		\\
		+
		\frac{1}{\varepsilon} \tD_{A_+}(\tu) \, \tMlumped \, \tD_{A_+}(\tu) \, (\tu - \tu_b)
		+
		\tF \, \tD(\ty,\tu)^{-2} \tMlumped \, \tP_\Omega \tp
		.
		\label{eq:first-order_derivative_Lagrangian_wrt_u_discrete}
	\end{multline}
\end{subequations}
\makeatother
Here $\tD_{A_+}(\tu)$ and $\tD_{A_-}(\tu)$ are diagonal (active-set) matrices with entries
\makeatletter
\begin{equation*}
	\begin{aligned}
		[\tD_{A_+}(\tu)]_{ii} 
		&
		= 
		\begin{cases}
			1 & \text{where } [\tu_a - \tu]_i \ge 0
			,
			\\
			0 & \text{otherwise}
			,
		\end{cases}
		&
		\qquad
		\ltx@ifclassloaded{svjour3}{\\}{}
		[\tD_{A_-}(\tu)]_{ii} 
		&
		= 
		\begin{cases}
			1 & \text{where } [\tu - \tu_b]_i \ge 0
			,
			\\
			0 & \text{otherwise}
			,
		\end{cases}
		\end{aligned}
\end{equation*}
\makeatother
and we set $\tD_A(\tu) = \tD_{A_+}(\tu) + \tD_{A_-}(\tu)$.

In order to solve the discrete optimality system consisting of \eqref{eq:state_equation_discrete_with_vectors} and \eqref{eq:first-order_derivatives_Lagrangian_discrete}, we employ a finite-dimensional semismooth Newton method (\cref{algorithm:semismooth_Newton_method_discrete}).
This requires the evaluation of first-order derivatives of the state equation \eqref{eq:state_equation_discrete_with_vectors} as well as second-order derivatives of the Lagrangian \eqref{eq:Lagrangian_discrete}.
The following expressions are obtained.
\makeatletter
\begin{subequations}
	\label{eq:blocks_of_Newton_matrix_discrete}
	\begin{align}
		e_\ty(\ty,\tu)
		&
		=
		\tP_\Omega \tK 
		+
		2 \, \tP_\Omega \, \tMlumped \, \tF \, \tB \, \tD(\ty,\tu)^{-2} \bone \, \ty^\transp \, \tK 
		+ 
		\tP_\Gamma 
		,
		\label{eq:first-order_derivative_state_equation_wrt_y_discrete}
		\\
		e_\tu(\ty,\tu)
		&
		=
		\tP_\Omega \tMlumped \, \tD(\ty,\tu)^{-2} \tF 
		,
		\label{eq:first-order_derivative_state_equation_wrt_u_discrete}
		\\
		\cL_{\ty\ty}(\ty,\tu,\tp)
		&
		=
		\tM
		- 8 \tp^\transp \tP_\Omega \tMlumped \, \tF \, \tB^2 \tD(\ty,\tu)^{-3} \bone \, \tK \, \ty \, \ty^\transp \tK
		\ltx@ifclassloaded{svjour3}{%
			\notag
			\\
			&
			\quad
		}{}
		+ 2 \tp^\transp \tP_\Omega \tMlumped \, \tF \, \tB \, \tD(\ty,\tu)^{-2} \bone \, \tK
		,
		\label{eq:second-order_derivative_Lagrangian_wrt_yy_discrete}
		\\
		\cL_{\ty\tu}(\ty,\tu,\tp)
		&
		=
		- 4 \, \tK \, \ty \, \tp^\transp \tP_\Omega \tMlumped \, \tF \, \tD(\ty,\tu)^{-3} \tB
		,
		\label{eq:second-order_derivative_Lagrangian_wrt_yu_discrete}
		\\
		\cL_{\tu\tu}(\ty,\tu,\tp)
		&
		=
		\lambda \, (\tK + \tM) 
		+
		\frac{1}{\varepsilon} \tD_A(\tu) \, \tMlumped \, \tD_A(\tu) 
		\ltx@ifclassloaded{svjour3}{%
			\notag
			\\
			&
			\quad
		}{}
		-
		2 \, \diag(\tMlumped \, \tp) \, \tD(\ty,\tu)^{-3} \tF
		.
		\label{eq:second-order_derivative_Lagrangian_wrt_uu_discrete}
	\end{align}
\end{subequations}
\makeatother
Notice that the expression for $\cL_{\tu\tu}$ is the generalized derivative of $\cL_\tu$ in the sense of \cref{definition:Newton_differentiability}.

The discrete generalized Newton system has the following form:
\begin{equation}
	\label{eq:generalized_Newton_sytem_discrete}
	\begin{bmatrix}
		\cL_{\ty\ty}(\ty,\tu,\tp) & \cL_{\ty\tu}(\ty,\tu,\tp) & e_\ty(\ty,\tu)^\transp 
		\\
		\cL_{\tu\ty}(\ty,\tu,\tp) & \cL_{\tu\tu}(\ty,\tu,\tp) & e_\tu(\ty,\tu)^\transp 
		\\
		e_\ty(\ty,\tu) & e_\tu(\ty,\tu) & \bnull
	\end{bmatrix}
	\begin{pmatrix}
		\delta \ty \\ \delta \tu \\ \delta \tp
	\end{pmatrix}
	=
	-
	\begin{pmatrix}
		\cL_\ty(\ty,\tu,\tp)
		\\
		\cL_\tu(\ty,\tu,\tp)
		\\
		e(\ty,\tu)
	\end{pmatrix}
	.
\end{equation}
In contrast to standard optimal control problems which do not feature a nonlocal PDE, some of the blocks in \eqref{eq:generalized_Newton_sytem_discrete} are no longer sparse.
This comment applies to $e_\ty$ due to the second summand in \eqref{eq:first-order_derivative_state_equation_wrt_y_discrete}, to $\cL_{\ty\ty}$ due to the second summand in \eqref{eq:second-order_derivative_Lagrangian_wrt_yy_discrete} as well as to $\cL_{\ty\tu}$ given by \eqref{eq:second-order_derivative_Lagrangian_wrt_yu_discrete}.
For a high performance implementation, it is therefore important to not assemble the blocks in \eqref{eq:generalized_Newton_sytem_discrete} as matrices, but rather to provide matrix-vector products and use a preconditioned iterative solver such as \minres (\cite{PaigeSaunders:1975:1}) to solve \eqref{eq:generalized_Newton_sytem_discrete}.
This aspect, however, is beyond the scope of this paper and we defer the design and analysis of a suitable preconditioner to future work.
For the time being we resort to the direct solution of \eqref{eq:generalized_Newton_sytem_discrete} using \matlab's direct solver, which is still feasible on moderately fine discretizations of two-dimensional domains.

Our implementation of the semismooth Newton method is described in \cref{algorithm:semismooth_Newton_method_discrete}. 
In contrast to \cref{algorithm:semismooth_Newton_method}, we added an additional step in which we solve the discrete nonlinear state equation \eqref{eq:state_equation_discrete} for $\ty_{k+1}$ once per iteration for increased robustness; see \cref{step:nonlinear_state_update} in \cref{algorithm:semismooth_Newton_method_discrete}.
Notice that the preliminary linear update to $\ty_{k+1}$ in \cref{step:linear_update_of_iterates} is still useful since it provides an initial guess for the subsequent solution of $e(\ty_{k+1},\tu_{k+1}) = 0$.
We mention that nonlinear state updates have been analyzed in the closely related context of SQP methods, \eg, in  \cite{Ulbrich:2007:1,CleverLangUlbrichZiems:2011:1}.
We also added a rudimentary damping strategy which improves the convergence behavior.
In our examples, it suffices to choose $\damping = 1/2$ when $\norm{\cL_\tu(\ty_k,\tu_k,\tp_k)}_{(\tK+\tM)^{-1}} > 1/10$ and $\damping = 1$ otherwise.

The stopping criterion we employ in \cref{step:stopping_criterion} measures the three components of the residual, \ie, the right-hand side in \eqref{eq:generalized_Newton_sytem_discrete}.
Following the function space setting of the continuous problem, we evaluate the (squared) $H^{-1}(\Omega)$-norm of all residual components, which amounts to 
\begin{equation}
	\label{eq:residual_norms}
	R^2(\ty,\tu,\tp)
	\coloneqq
	\norm{\cL_\ty(\ty,\tu,\tp)}_{(\tK+\tM)^{-1}}^2
	+
	\norm{\cL_\tu(\ty,\tu,\tp)}_{(\tK+\tM)^{-1}}^2
	+
	\norm{e(\ty,\tu)}_{(\tK+\tM)^{-1}}^2
	.
\end{equation}
\Cref{algorithm:semismooth_Newton_method_discrete} is stopped when
\begin{equation}
	\label{eq:stopping_criterion}
	R(\ty,\tu,\tp)
	\le
	10^{-6}
\end{equation}
is reached.
Moreover, we impose a tolerance of $\norm{e(\ty,\tu)}_{(\tK+\tM)^{-1}} \le 10^{-10}$ for the solution of the forward problem in \cref{step:nonlinear_state_update}.

\begin{algorithm}[Discrete semismooth Newton method with nonlinear state update for the solution of a discretized instance of problem \eqref{eq:optimal_control_problem_penalized}] \hfill
	\label{algorithm:semismooth_Newton_method_discrete}
	\begin{algorithmic}[1]
		\Require initial guess $(\ty_0,\tu_0,\tp_0) \in V_h \times V_h \times V_h$ 
		\Ensure approximate stationary point of the discretized instance of \eqref{eq:optimal_control_problem_penalized}
		\State Set $k \coloneqq 0$
		\While{not converged}
		\label{step:stopping_criterion}
		\State Determine the active sets $A_+(\tu_k)$ and $A_-(\tu_k)$
		\State Solve the Newton system \eqref{eq:generalized_Newton_sytem_discrete} for $(\delta \ty, \delta \tu, \delta \tp)$, given $(\ty_k,\tu_k,\tp_k)$
		\label{step:solve_SNN_system}
		\State Update the iterates by setting
		\begin{equation*}
			\ty_{k+1} 
			\coloneqq
			\ty_k + \damping \, \delta \ty
			,
			\quad
			\tu_{k+1} 
			\coloneqq
			\tu_k + \damping \, \delta \tu
			,
			\quad
			\tp_{k+1} 
			\coloneqq
			\tp_k + \damping \, \delta \tp
		\end{equation*}
		where $\damping \in (0,1]$ is a suitable damping parameter.
		\label{step:linear_update_of_iterates}
		\State Solve the nonlinear state equation \eqref{eq:state_equation_discrete} for the state $\ty_{k+1}$, given the control $\tu_{k+1}$ 
		\label{step:nonlinear_state_update}
		\State Set $k \coloneqq k+1$
		\EndWhile
	\end{algorithmic}
\end{algorithm}

\section{Numerical Experiments}
\label{section:numerical_experiments}

In this section we describe a number of numerical experiments.
The first experiment serves the purpose of demonstrating the influence of the non-locality parameter~$b$.
In the second experiment, we numerically confirm the mesh independence of our algorithm. 
The third experiment is dedicated to studying the impact of the penalty parameter~$\varepsilon$.

In a slight extension of \eqref{eq:objective_discrete}, we distinguish two control cost parameters $\lambda_1$ and $\lambda_2$, which leads to discrete problems of the form
\begin{align}
	\label{eq:objective_discrete_with_two_control_cost_parameters}
	J(\ty,\tu)
	&
	=
	\frac{1}{2} (\ty - \ty_d)^\transp \tM (\ty - \ty_d)
	+
	\frac{\lambda_1}{2} \tu^\transp \tK \tu
	+
	\frac{\lambda_2}{2} \tu^\transp \tM \tu
	\notag
	\\
	& 
	\quad
	+
	\frac{1}{2 \varepsilon} (\tu_a - \tu)_+^\transp \tMlumped (\tu_a - \tu)_+
	+
	\frac{1}{2 \varepsilon} (\tu - \tu_b)_+^\transp \tMlumped (\tu - \tu_b)_+
	.
\end{align}
The modifications to \eqref{eq:Lagrangian_discrete}--\eqref{eq:blocks_of_Newton_matrix_discrete} due to the two control cost parameters are obvious.
As mentioned in \cref{section:discretization}, our implementation of \cref{algorithm:semismooth_Newton_method_discrete} employs a direct solver for the linear systems arising in \cref{step:solve_SNN_system} and is therefore only suitable for relatively coarse discretization of two-dimensional domains.
Unless otherwise mentioned, the following experiments are obtained on a mesh discretizing a square domain with $N_V = 665$~vertices and $N_T = 1248$~triangles.
Notice that convex domains are covered by our theory due to \cref{remark:convex_domains}.
The typical run-time for \cref{algorithm:semismooth_Newton_method_discrete} is around \SI{3}{\second}.

\subsection{Influence of the Non-Locality Parameter}
\label{subsection:numerical_experiments_non-locality}

Our initial example builds on the two-dimensional problem presented in \cite{DelgadoFigueiredoGayteMoralesRodrigo:2017:1}.
The problem domain is $\Omega = (-0.5,0.5)^2$; notice that this is slightly incorrectly stated in \cite{DelgadoFigueiredoGayteMoralesRodrigo:2017:1}. 
Moreover, we have right-hand side $f(x,y) \equiv 100$ and desired state $y_d(x,y) \equiv 0$.
The lower bound for the control is given as $u_a(x,y) = -3x -3y +10$ and the upper bound is $u_b \equiv \infty$.
Moreover, the control cost parameters are $\lambda_1 = 0$ and $\lambda_2 = 4 \cdot 10^{-5}$.
We choose $\varepsilon = 10^{-2}$ as our penalty parameter.
The coefficient function determining the degree of non-locality is set to $b(x,y) = \alpha \, (x^2 + y^2)$, where $\alpha$ varies in $\{0,10^0,10^1,10^2,10^3\}$.
The case $\alpha = 1$ is considered in \cite{DelgadoFigueiredoGayteMoralesRodrigo:2017:1} and we reproduce their results.

For each value of $\alpha$, we start from an initial guess constructed as follows.
We initialize $\tu_0$ to the lower bound~$\tu_a$ and set $\ty_0$ to the numerical solution of the forward problem with control $\tu_0$. 
The adjoint state is initialized to $\tp_0 = \bnull$.

\begin{figure}[htp]
	\begin{subfigure}{0.49\textwidth}
		\includegraphics[width=\linewidth]{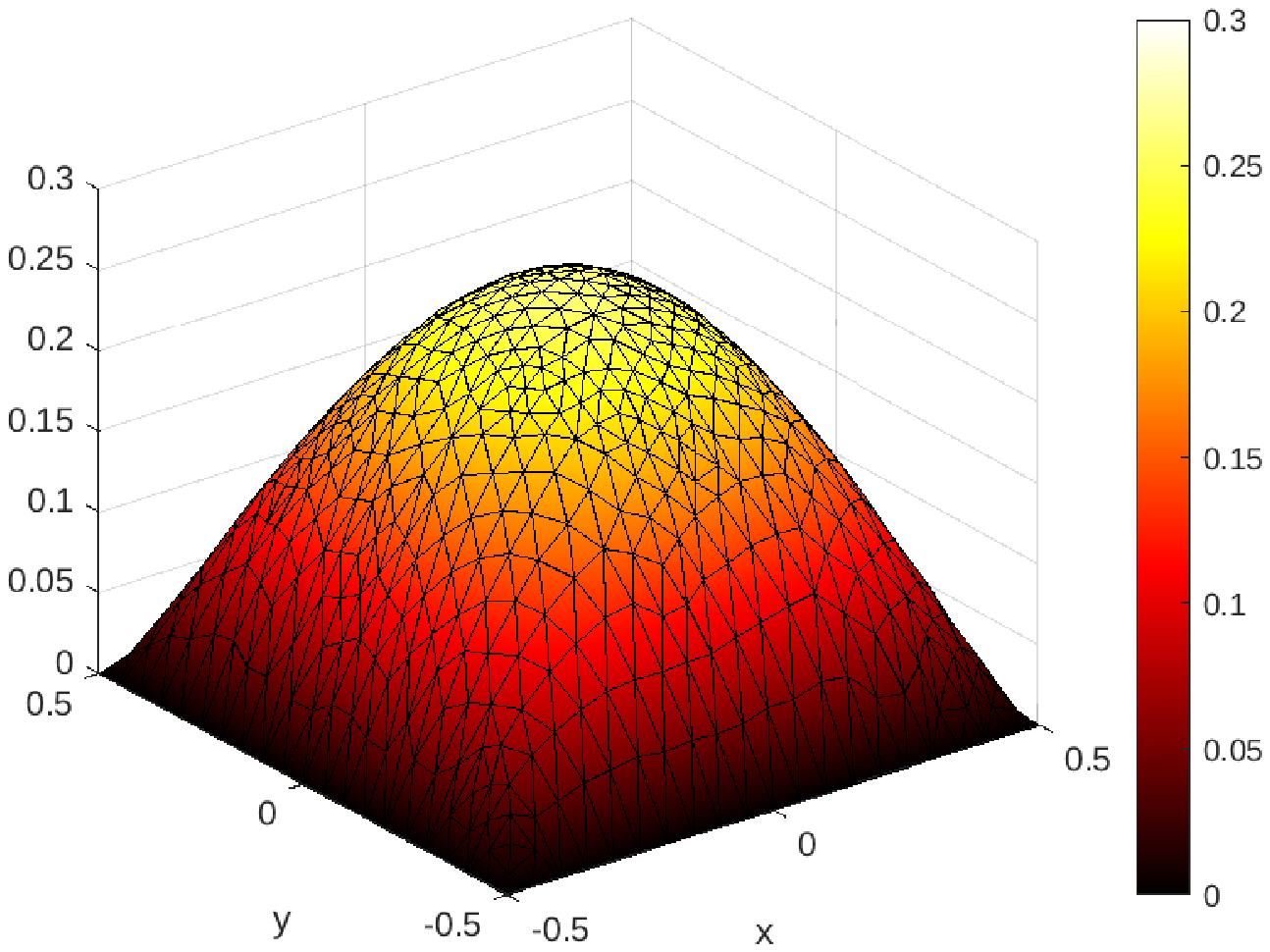}
	\end{subfigure}
	\hfill
	\begin{subfigure}{0.49\textwidth}
		\includegraphics[width=\linewidth]{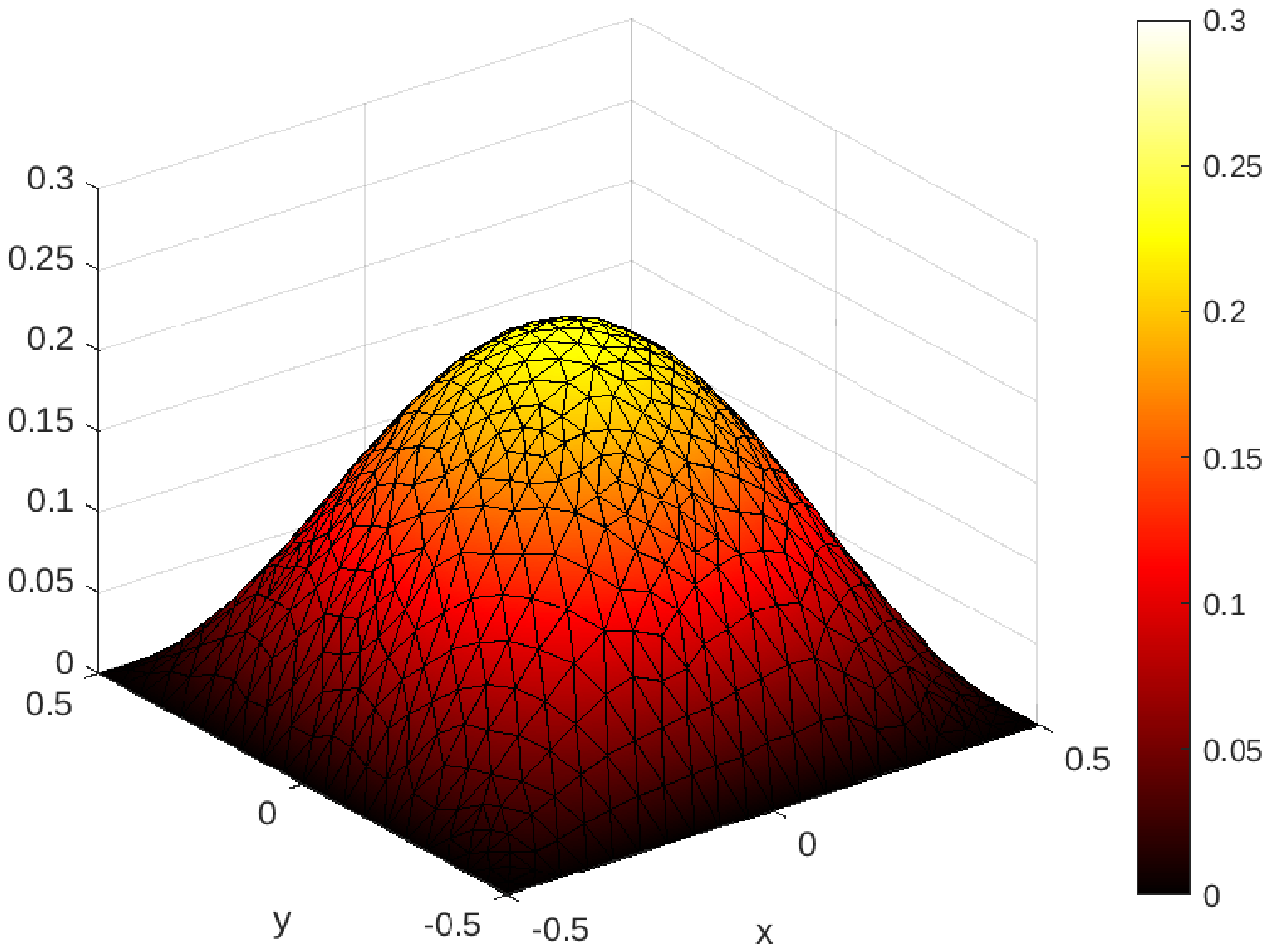}
	\end{subfigure}
	\begin{subfigure}{0.49\textwidth}
		\includegraphics[width=\linewidth]{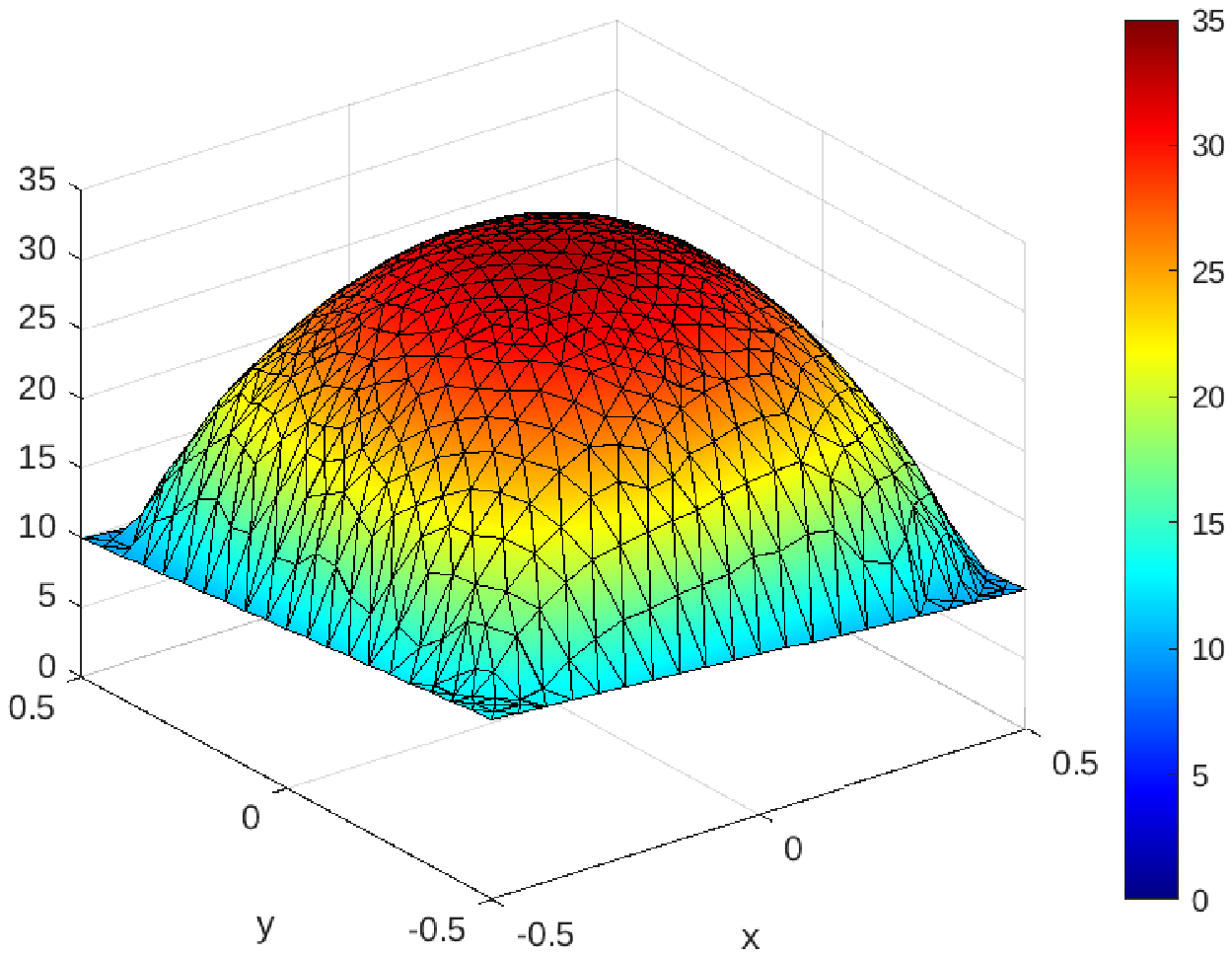}
	\end{subfigure}
	\hfill
	\begin{subfigure}{0.49\textwidth}
		\includegraphics[width=\linewidth]{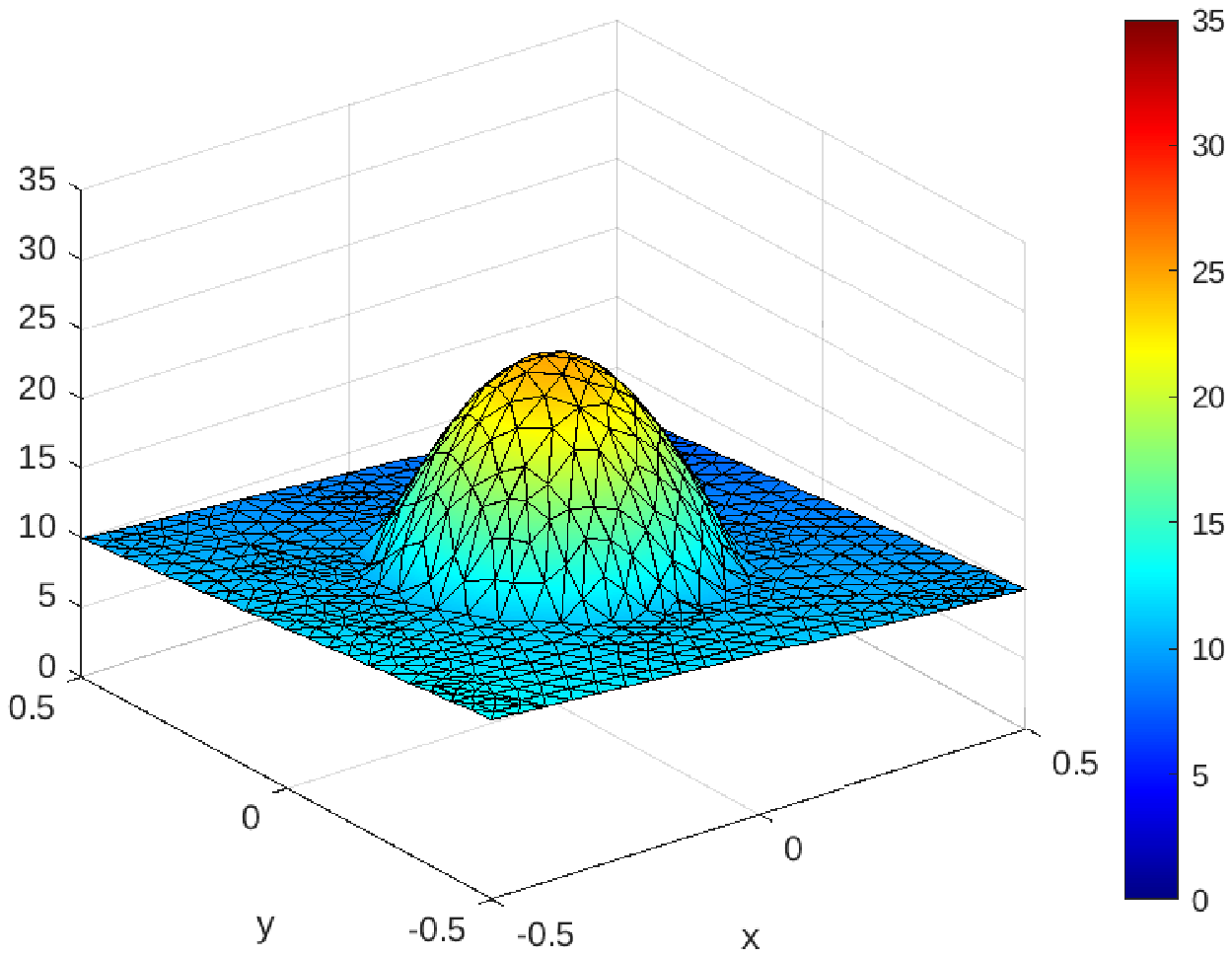}
	\end{subfigure}
	\begin{subfigure}{0.49\textwidth}
		\includegraphics[width=0.85\linewidth]{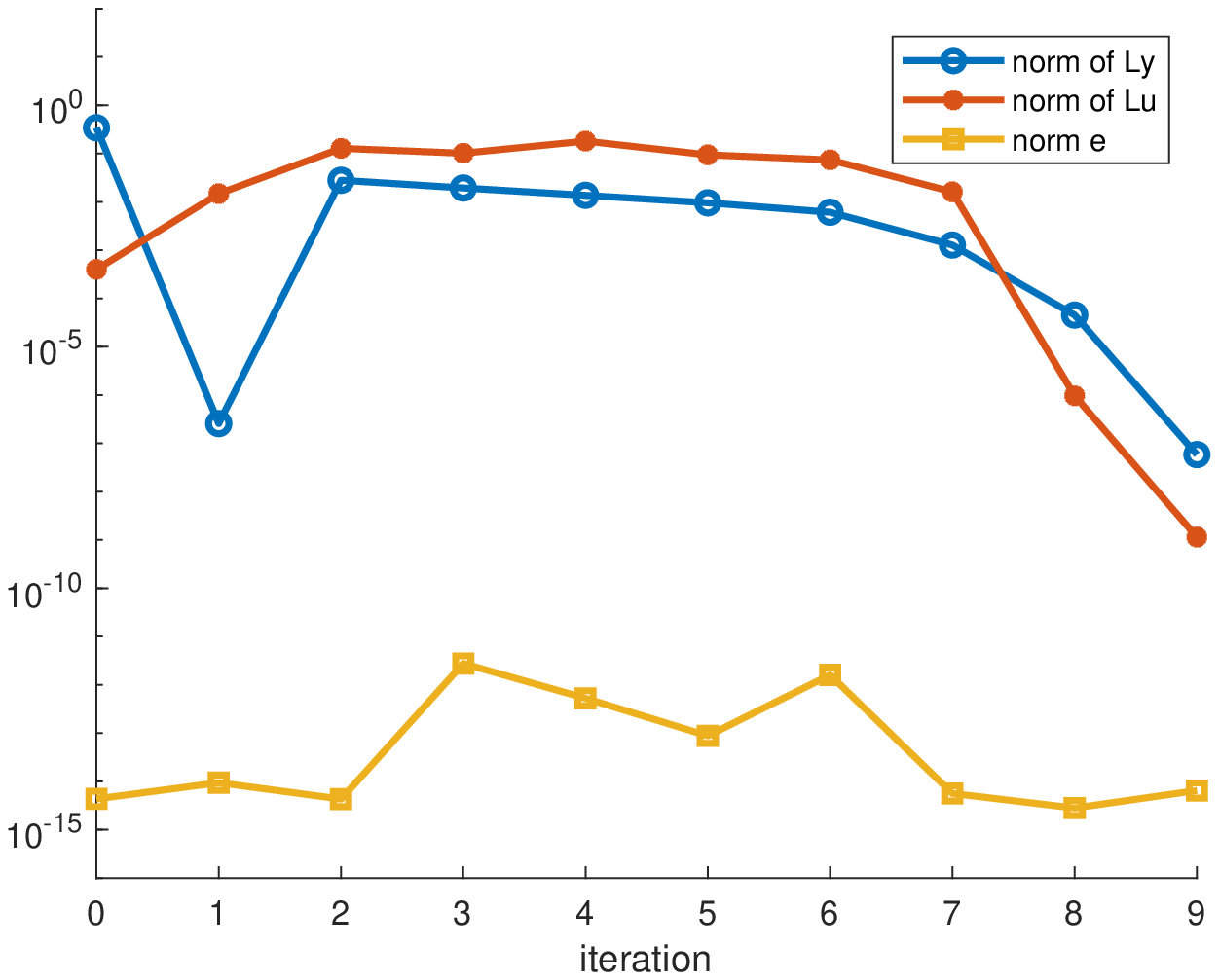}
	\end{subfigure}
	\hfill
	\begin{subfigure}{0.49\textwidth}
		\includegraphics[width=0.85\linewidth]{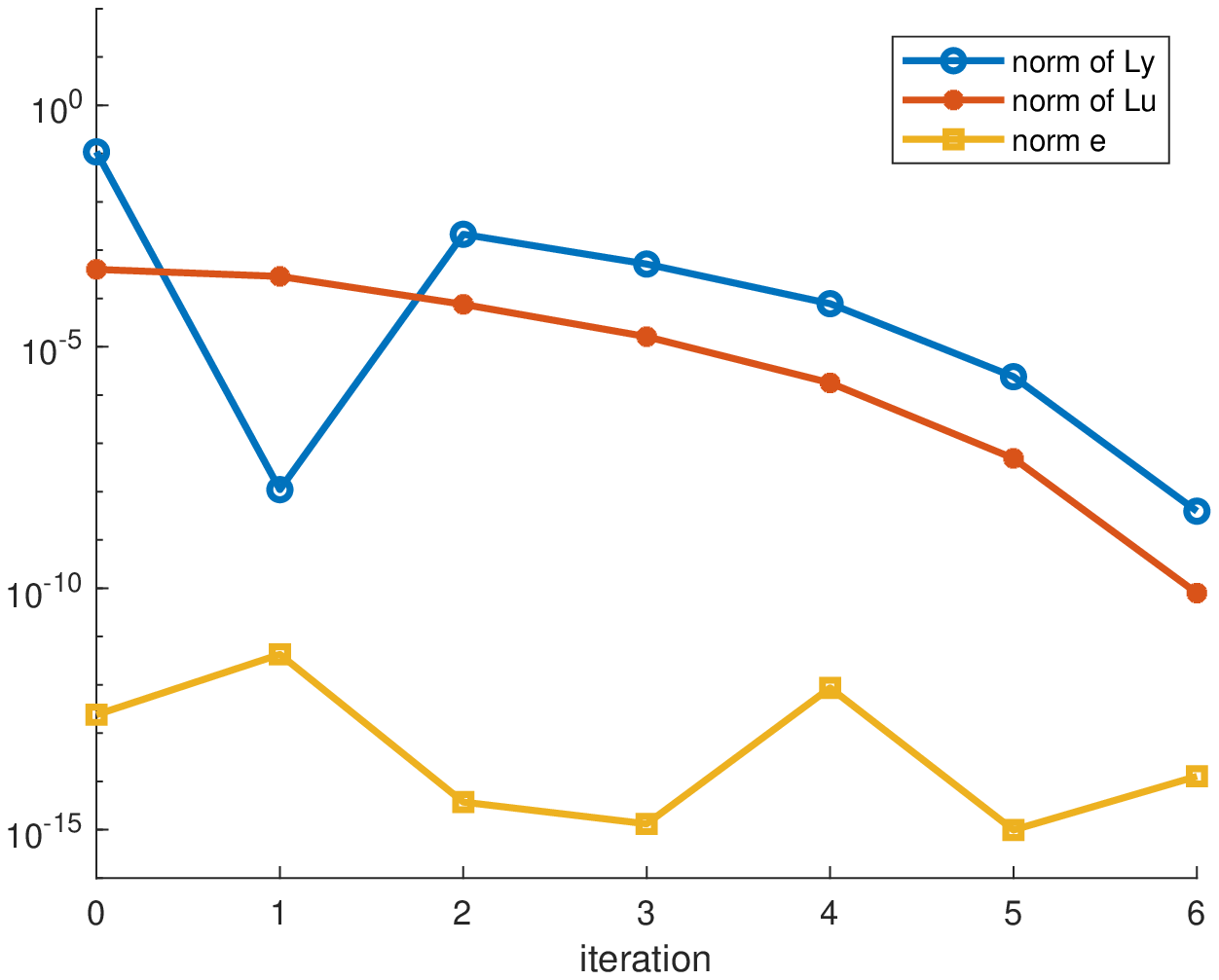}
	\end{subfigure}
	\caption{Optimal states $y$ (top row), optimal controls $u$ (middle row) and convergence history (bottom row) obtained for the example from \cref{subsection:numerical_experiments_non-locality} for $\alpha = 1$ (left column) and $\alpha = 10^3$ (right column). The three norms shown in the convergence plots correspond to the three terms in \eqref{eq:residual_norms}, \ie, $\norm{\cL_\ty(\ty,\tu,\tp)}_{(\tK+\tM)^{-1}}$, $\norm{\cL_\tu(\ty,\tu,\tp)}_{(\tK+\tM)^{-1}}$ and $\norm{e(\ty,\tu)}_{(\tK+\tM)^{-1}}$.}
	\label{figure:numerical_experiments_non-locality_state_control_history}
\end{figure}

\Cref{figure:numerical_experiments_non-locality_state_control_history} show some of the optimal state and control functions obtained.
We notice that the solution in case of a local problem ($\alpha = 0$) is visually indistinguishable from the setting $\alpha = 1$ considered in \cite{DelgadoFigueiredoGayteMoralesRodrigo:2017:1}.
We therefore compare it to the case $\alpha = 10^3$ of significantly more pronounced non-local effects.
Clearly, an increase in the non-local parameter aids the control in this example, so the control effort can decrease, as reflected in \cref{figure:numerical_experiments_non-locality_state_control_history}.
Also, we observe that the number of iterations of the discrete semismooth Newton method (\cref{algorithm:semismooth_Newton_method_discrete}) decreases slightly as $\alpha$ increases; see \cref{table:numerical_experiments_non-locality_number_of_SSN_steps}.

\begin{table}[ht]
	\centering
	\begin{tabular}{@{}rr@{}}
  \toprule
           $\alpha$ &          iterations \\
  \midrule
  0.00e+00 &         10 \\
  1.00e+00 &          9 \\
  1.00e+01 &          7 \\
  1.00e+02 &          7 \\
  1.00e+03 &          6 \\
  \bottomrule
\end{tabular}

	\caption{Number of iterations of the discrete semismooth Newton method (\cref{algorithm:semismooth_Newton_method_discrete}) for various values of the non-locality parameter~$\alpha$ in the example from \cref{subsection:numerical_experiments_non-locality}.}
	\label{table:numerical_experiments_non-locality_number_of_SSN_steps}
\end{table}

\subsection{Dependence on the Discretization}
\label{subsection:numerical_experiments_mesh_size}

In this experiment we study the dependence of the number of semismooth Newton steps in \cref{algorithm:semismooth_Newton_method_discrete} on the refinement level of the underlying discretization.
To this end, we consider a coarse mesh and two uniform refinements; see \cref{table:numerical_experiments_mesh_size_number_of_SSN_steps}.

The problem is similar as in \cref{subsection:numerical_experiments_non-locality}.
The domain is $\Omega = (-0.5,0.5)^2$.
We use $f(x,y) \equiv 100$ as right-hand side and the desired state is $y_d(x,y) \equiv 0$.
The lower bound for the control is now given as $u_a(x,y) = -10x -10y +20$ and the upper bound is $u_b = u_a + 5$.
Moreover, the control cost parameters are $\lambda_1 = 0$ and $\lambda_2 = 4 \cdot 10^{-5}$.
We choose $\varepsilon = 10^{-2}$ as our penalty parameter.
The coefficient function determining the degree of non-locality is set to $b(x,y) = 100 \, (x^2 + y^2)$.

\begin{figure}[htp]
	\centering
	\begin{subfigure}{0.45\textwidth}
		\includegraphics[width=\linewidth]{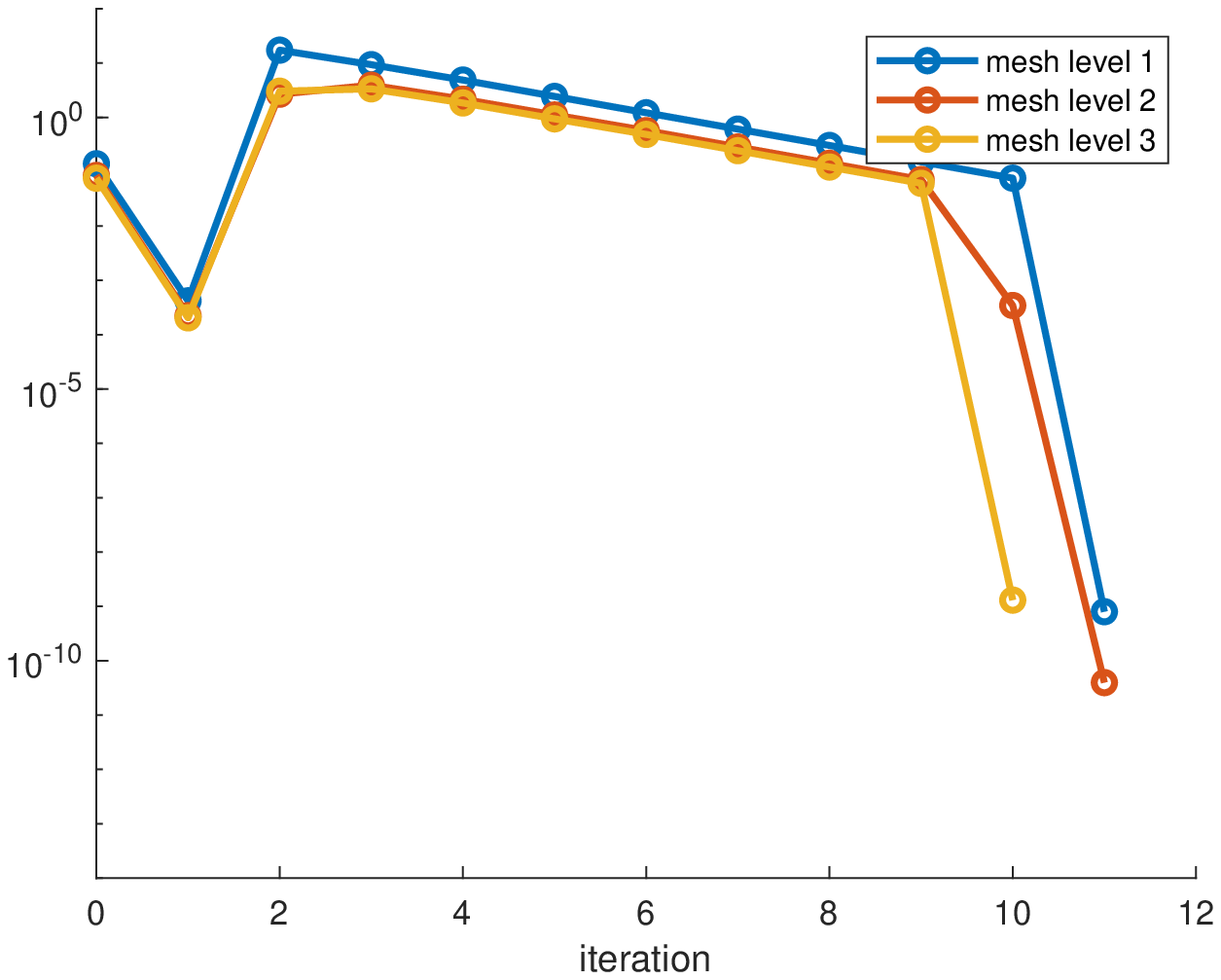}
	\end{subfigure}
	\begin{subfigure}{0.49\textwidth}
		\includegraphics[width=\linewidth]{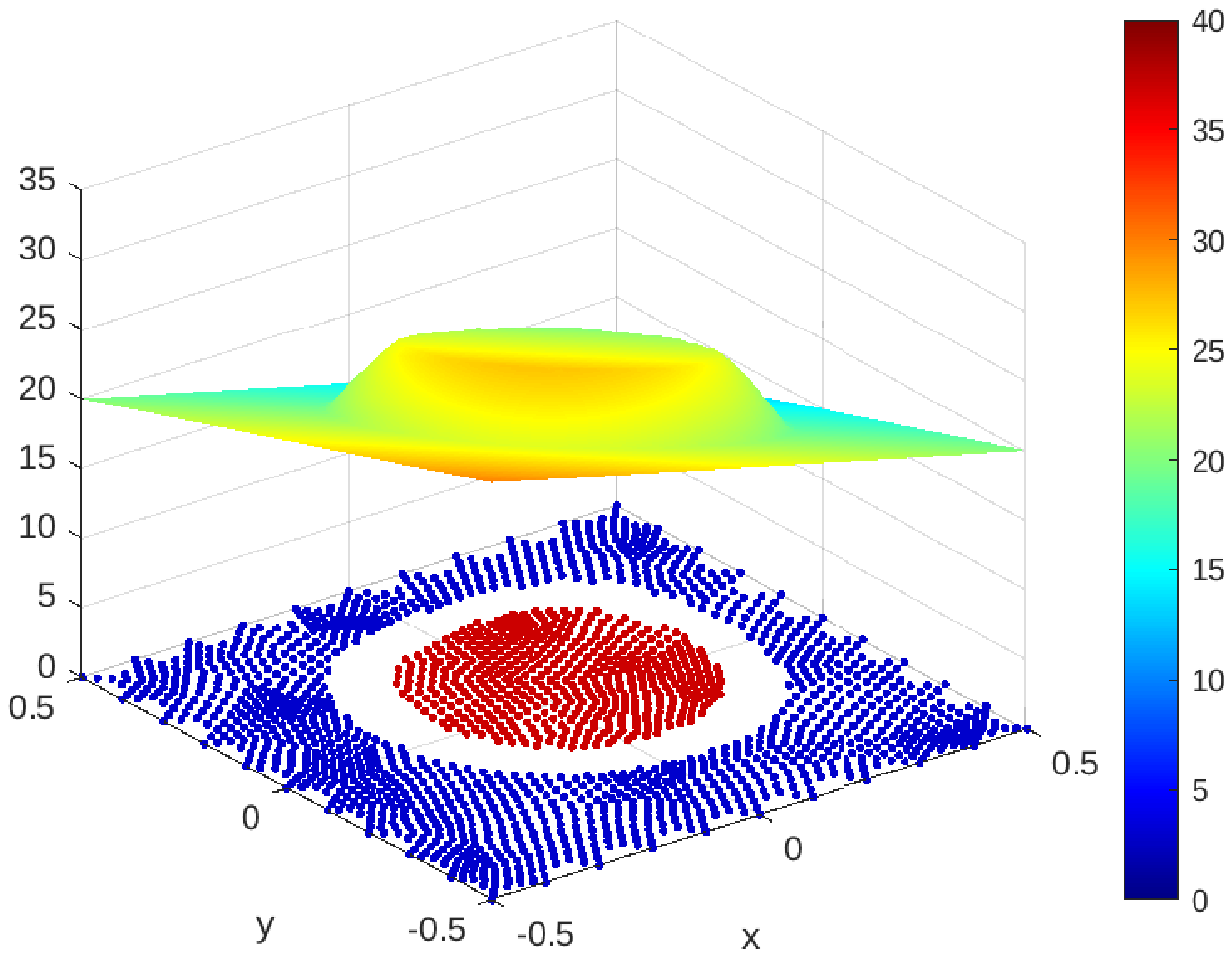}
	\end{subfigure}
	\caption{The convergence plot (left column) shows the total residual norm $R(\ty,\tu,\tp)$ as in \eqref{eq:residual_norms} on all mesh levels for the example from \cref{subsection:numerical_experiments_mesh_size}. The control on the finest level is shown in the right column. Nodes where $u = u_b$ and $u = u_a$ holds are shown in red and blue, respectively.}
	\label{figure:numerical_experiments_mesh_size}
\end{figure}

For each mesh, we start from an initial guess constructed as follows.
We initialize $\tu_0$ to the lower bound~$\tu_a$ and set $\ty_0$ to the numerical solution of the forward problem with control $\tu_0$. 
The adjoint state is initialized to $\tp_0 = \bnull$.
In this example, both the lower and upper bounds are relevant on all mesh levels.
Nonetheless, we observe a mesh-independent convergence behavior; see \cref{figure:numerical_experiments_mesh_size}.

\begin{table}[ht]
	\centering
	\begin{tabular}{@{}rrrr@{}}
  \toprule
           level &          $N_V$ &          $N_T$ &          iterations \\
  \midrule
      1 &   177 &   312 &         11 \\
      2 &   665 &  1248 &         11 \\
      3 &  2577 &  4992 &         10 \\
  \bottomrule
\end{tabular}

	\caption{Number of iterations of the discrete semismooth Newton method (\cref{algorithm:semismooth_Newton_method_discrete}) for various mesh levels in the example from \cref{subsection:numerical_experiments_mesh_size}.}
	\label{table:numerical_experiments_mesh_size_number_of_SSN_steps}
\end{table}

\subsection{Influence of the Penalty Parameters}
\label{subsection:numerical_experiments_penalty}

In this final experiment, we study the behavior of \cref{algorithm:semismooth_Newton_method_discrete} and the solutions to the penalized problem \eqref{eq:optimal_control_problem_penalized} in dependence of the penalty parameter~$\varepsilon$.
	We solve similar problems as before, with domain $\Omega = (-0.5,0.5)^2$, right-hand side $f(x,y) \equiv 100$ and desired state $y_d(x,y) \equiv 0$.
	The lower bound for the control is $u_a(x,y) = -10x -10y +20$ and the upper bound is $u_b = u_a + 8$.
Moreover, the control cost parameters are $\lambda_1 = 0$ and $\lambda_2 = 4 \cdot 10^{-5}$.
The penalty parameter varies in $\{10^0, 10^{-1}, 10^{-2}, 10^{-3}, 10^{-4}\}$.
The coefficient function determining the degree of non-locality is set to $b(x,y) = 100 \, (x^2 + y^2)$.

\begin{figure}[htp]
	\centering
	\begin{subfigure}{0.49\textwidth}
		\includegraphics[width=\linewidth]{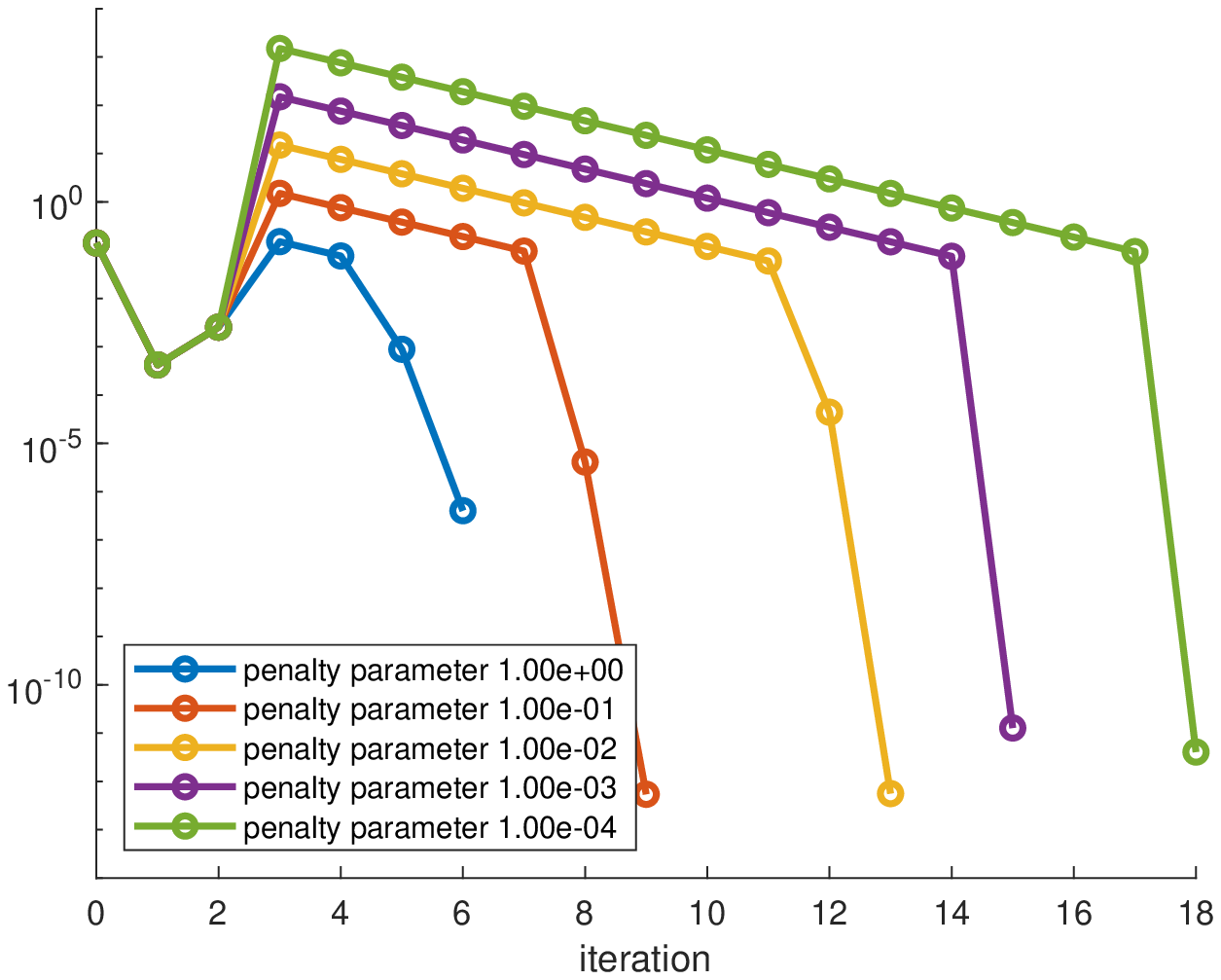}
	\end{subfigure}
	\begin{subfigure}{0.49\textwidth}
		\includegraphics[width=\linewidth]{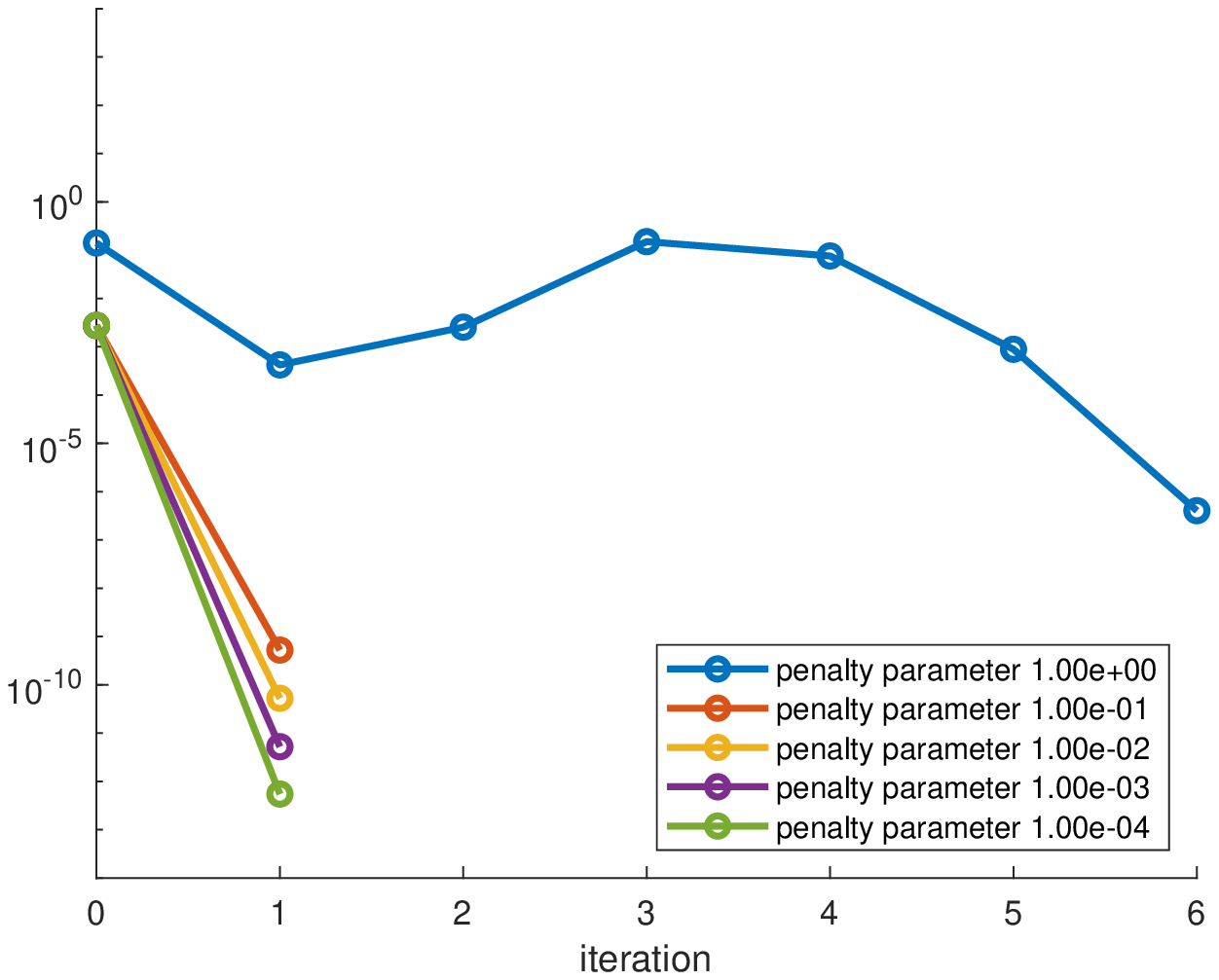}
	\end{subfigure}
	\caption{The convergence plot shows the total residual norm $R(\ty,\tu,\tp)$ as in \eqref{eq:residual_norms} for all values of the penalty parameter~$\varepsilon$. In the left plot, the same initial guess was used for all penalty parameters. With warmstarting, convergence can be achieved in one semismooth Newton step.}
	\label{figure:numerical_experiments_penalty}
\end{figure}

The construction of an initial guess is the same as in \cref{subsection:numerical_experiments_mesh_size}.
The experiment is split into two parts.
First, we consider \cref{algorithm:semismooth_Newton_method_discrete} without warmstarts.
The corresponding results are shown in \cref{table:numerical_experiments_penalty_number_of_SSN_steps}.
As expected, the number of Newton steps increases as $\varepsilon \searrow 0$ while the norm of the bound violation decreases.
Second, we repeat the same experiment with warmstarts.
That is, we use the initialization as described above only for the initial value of $\varepsilon$.
Subsequent runs of \cref{algorithm:semismooth_Newton_method_discrete} are initialized with the final iterates obtained for the previous value of $\varepsilon$.
This strategy is very effective, as shown in \cref{figure:numerical_experiments_penalty} (right column).

\begin{table}[ht]
	\centering
	\begin{tabular}{@{}rrrr@{}}
  \toprule
      $\varepsilon$ &          iterations & $\norm{(u_a - u)_+}_{L^\infty(\Omega)}$ & $\norm{(u - u_b)_+}_{L^\infty(\Omega)}$ \\
  \midrule
  1.00e+00 &          6 &                       1.19e-03 &                       4.41e-04 \\
  1.00e-01 &          9 &                       1.19e-04 &                       4.41e-05 \\
  1.00e-02 &         13 &                       1.19e-05 &                       4.41e-06 \\
  1.00e-03 &         15 &                       1.19e-06 &                       4.41e-07 \\
  1.00e-04 &         18 &                       1.19e-07 &                       4.41e-08 \\
  \bottomrule
\end{tabular}

	\caption{Number of iterations of the discrete semismooth Newton method (\cref{algorithm:semismooth_Newton_method_discrete}, without warmstart) for various values of the penalty parameter~$\varepsilon$ in the example from \cref{subsection:numerical_experiments_penalty}. The terms $\norm{(u_a - u)_+}_{L^\infty(\Omega)}$ and $\norm{(u - u_b)_+}_{L^\infty(\Omega)}$ refer to the maximal positive nodal values of $u_a - u$ and $u - u_b$, respectively.}
	\label{table:numerical_experiments_penalty_number_of_SSN_steps}
\end{table}

\appendix
\section{Comment on the Proof of Existence of an Optimal Solution}
\label{section:appendix}

We believe that the proof concerning the existence of an optimal solution in Theorem~2.5 of \cite{DelgadoFigueiredoGayteMoralesRodrigo:2017:1} contains a flaw.
Indeed, step~4 of the proof implies that for every weakly convergent sequence $\{a_n\}$ in $L^2(\Omega)$ and for each $\varepsilon > 0$, there exists a set $A_\varepsilon$ with $\abs{A_\varepsilon} < \varepsilon$ and a subsequence which converges pointwise in a dense subset of $\Omega\setminus A_\varepsilon$.
The following counterexample shows that this is not the case.

\begin{example*}
	Let $f(x)$ be a $1$-periodic function on $\R$ and
	\begin{equation*}
		f(x)
		\coloneqq
		\begin{cases}
			0, & 0 \le x \le 1/2 
			,
			\\
			1, & 1/2 < x \le 1
			.
		\end{cases}
	\end{equation*}  
	Set $a_n(x) = f(nx)$. 
	We can see
	\begin{equation*}
		a_n 
		\weakly 
		\int_0^1 f(x) \dx 
		=
		1/2
		,
	\end{equation*}
	see \cite[Theorem~2.6]{CioranescuDonato:1999:1}, but we cannot find a subsequence of $\{a_n\}$ which converges pointwise to $1/2$ for any $x \in [0,1]$.
\end{example*}

\section*{Acknowledgments}

MH would like to thank Morteza Fotouhi (Sharif University of Technology) for fruitful discussions concerning the material in \cref{section:optimality_system}.

\printbibliography

\end{document}